\documentclass[article,reqno]
{amsart}
\usepackage[utf8]{inputenc}
\usepackage{a4wide, amsfonts}
\usepackage[english]{babel}
\usepackage{xcolor}
\usepackage[T1]{fontenc}
\usepackage{amsmath, amssymb,epic,graphicx,mathrsfs,enumerate}    
\usepackage{cite}
\usepackage[numbers,sort&compress]{natbib}
\usepackage[all]{xy}
\usepackage{color}
\usepackage{comment}
\usepackage{enumitem}
\usepackage{esint}
\usepackage{hyperref}
\usepackage{amsthm}
\usepackage{latexsym}
\usepackage{epsfig}
\usepackage{soul}
\usepackage{geometry}
\usepackage{doi}
\hypersetup{colorlinks=true, linkcolor=black, filecolor=black, urlcolor=blue,
citecolor=blue}

\newtheorem{theorem}{Theorem}[section]

\newtheorem{lemma}[theorem]{Lemma}
\newtheorem{prop}[theorem]{Proposition}
\newtheorem{definition}{Definition}[section]
\theoremstyle{remark}
\newtheorem{obs}{Remark}

\newcommand{\R}{\mathbb{R}}

\newcommand{\E}{\mathcal{E}}
\newcommand{\e}{\varepsilon}
\newcommand{\al}{\alpha}
\newcommand{\g}{\gamma}
\newcommand{\la}{\lambda}
\newcommand{\s}{\sigma}
\newcommand{\p}{\varphi}
\newcommand{\de}{\delta}

\title[Ergodic MFG of Choquard-type]{Ergodic Mean-Field Games with aggregation of Choquard-type}

\author[C.Bernardini - A.Cesaroni]{Chiara Bernardini and Annalisa Cesaroni}
\address{Dipartimento di Matematica "Tullio Levi Civita",
Università di Padova, Via Trieste 121, 35121 Padova, Italy}
\email{annalisa.cesaroni@unipd.it, chiara.bernardini@math.unipd.it}
 
\begin{document}
\begin{abstract}
We consider second-order ergodic Mean-Field Games systems in the whole space $\R^N$ with coercive potential and aggregating nonlocal coupling, defined in terms of a Riesz interaction kernel.
These MFG systems describe  Nash equilibria of games with a large population of indistinguishable rational players attracted toward regions where the population is highly distributed. 
Equilibria solve a system of PDEs where an Hamilton-Jacobi-Bellman equation is combined with a  Kolmogorov-Fokker-Planck equation for the mass distribution. 
Due to the interplay  between the strength of the attractive term and the behavior of the diffusive part, we will obtain  three different regimes for the existence 
and non existence of classical solutions to the  MFG system. By means of a Pohozaev-type identity, we prove nonexistence of regular solutions to the MFG system without potential in the Hardy-Littlewood-Sobolev-supercritical regime. On the other hand,  using a fixed point argument,  we show existence of classical solutions in the Hardy-Littlewood-Sobolev-subcritical regime at least  for masses smaller than a given threshold value. In the mass-subcritical regime we  show that actually this threshold can be taken to be $+\infty$.\\


\end{abstract}
\maketitle

\date{}



\section{Introduction}

In this paper, we study ergodic Mean-Field Games systems defined  in the whole  space $\R^N$ with a coercive potential $V$ and attractive nonlocal coupling, defined in terms of a Riesz interaction kernel. More in details, given $M>0$, we consider elliptic systems of the form
\vspace{2mm}
\begin{equation}\label{1}
\begin{cases}
-\Delta u+\frac{1}{\gamma}|\nabla u|^\gamma+\la=V(x)-\int_{\R^N} \frac{m(y)}{|x-y|^{N-\alpha} }dy \\
-\Delta m-\mathrm{div} (m\nabla u(x)\,|\nabla u(x)|^{\g-2} )=0\\
\int_{\R^N}m=M,\quad m\ge0
\end{cases}\text{in}\,\,\R^N
\end{equation}
where $\g>1, \al\in(0,N)$ are fixed. Note that the unknowns in the system \eqref{1} are  the functions $u,m$ and the constant $\la\in\R$ which can be interpreted as a  Lagrange multiplier, related to the mass constraint $\int_{\R^N}m=M$. 

We will assume that the potential $V$ is a locally H\"older continuous coercive function, that is there exist $b$ and $C_V$ positive constants such that 
\begin{equation}\label{V}
C_V^{-1}(\max\{|x|-C_V,0\})^b\le V(x)\le C_V(1+|x|)^b,\quad \forall x\in\R^N.
\end{equation}
 The assumption of $V$ to be non-negative is not restrictive, we can assume more generally that $V$ is bounded from below and shift appropriately $\lambda$. 
 
Finally, the coupling in the system is given through the interaction term $-m\ast K_\al$, where $K_\al$ is the Riesz potential of order $\al\in(0,N)$ defined as 
$$K_\al(x)=\frac{1}{|x|^{N-\al}}.$$
We assume the Hamiltonian in the system \eqref{1}  has  the form $H(p)= \frac{1}{\gamma}|p|^\gamma$  for sake of simplicity but actually it  may be more general, namely we may assume 
that   $H:\R^N\to\R$ is strictly convex, $H\in C^2(\R^N\setminus\{0\})$ and there exist $C_H,\,K>0$ and $\g>1$ such that $\forall p\in \R^N$, it holds 
$$C_H|p|^\g-K\le H(p)\le C_H|p|^\g$$
$$\nabla H(p)\cdot p-H(p)\ge K^{-1}|p|^\g-K$$
$$|\nabla H(p)|\le K |p|^{\g-1}. $$

Mean-Field Games have been introduced  in the seminal papers of Lasry and Lions \cite{LL1} and by Huang, Caines and Malhamé \cite{HCM} in order to describe Nash equilibria of differential games with infinitely many infinitesimal rational players; this led to a broader study, also encouraged by their powerful applications in a wide range of disciplines (equations of this kind arise in Economics, Finance and models of social systems). The key idea underlying the theory comes from Statistical Mechanics and Physics, and consists in a mean-field approach to describe equilibria in a system of many interacting particles. The theory of Mean-Field Games models the behavior of a very large number of rational  and indistinguishable   players aiming at minimizing a certain  cost, by anticipating  the distribution of the overall population which result from the actions of all other players. We refer to \cite{GLL,GoPibook} for a general presentations of Mean-Field Games and their applications. In our setting, the dynamics of each player is described by the following controlled stochastic differential equation
$$dX_t=-v_t\,dt+\sqrt{2}\,dB_t$$
where $v_t$ is the controlled velocity and $B_t$ is a standard $N$ dimensional  Brownian motion. Each agent chooses $v_t$ in order to minimize the following long time average  cost 
$$ \lim\limits_{T\to\infty}\frac{1}{T}\mathbb{E}\int_0^T\left[\frac{|v_t|^{\g'}}{\g'}+V(X_t)-K_\al\ast m(X_t)\right]dt$$
where $\g'=\frac{\g}{\g-1}$ is the conjugate exponent of $\g$ and $m(x)$ is the density of population at $x\in\R^N$. In the ergodic setting, the distribution law of each player moving with  optimal speed converges as $t\to +\infty$ to an invariant measure $\mu$ (independent of the initial position) and $\mu$   coincides, in a equilibrium regime for the game, with the density of the population $m$.   From a PDE viewpoint, equilibria of the differential game are encoded by solutions of the system \eqref{1}, where the    Hamilton-Jacobi-Bellman equation takes into account the value of the game $\la$ and the optimal speed  $-\nabla u |\nabla u|^{\gamma-2}$ of the optimal control problem of a typical agent, and the Kolmogorov-Fokker-Planck equation gives the density of the overall  population $m$.\\

In the case when $\g=\g'=2$, as pointed out in \cite{LL1}, using the Hopf-Cole transformation $v(x):= e^{-u(x)/2}$ we can reduce the MFG system \eqref{1} to a single PDE. In particular we observe that with the previous change of variable, setting $m(x)=v^2(x)$, the MFG system \eqref{1} is equivalent to the normalized Choquard equation 
\begin{equation}\label{cho} \begin{cases}
-2\Delta v+(V(x)-\la)v=(K_\al\ast v^2)v\\
\int_{\R^N}v^2(x) dx=M, \qquad v>0
\end{cases}\quad\text{in}\,\,\R^N,\end{equation} 
with associated energy  
$$\E(v)=\int_{\R^N}2 |\nabla v|^2+ V(x)v^2 dx-\frac{1}{2}\int_{\R^N}\int_{\R^N} \frac{v^2 (x)\,v^2(y)} {|x-y|^{N-\al}}dx\,dy.$$
Choquard-type equations have been intensively studied during the last decades and have appeared in the context of various mean-field type physical models (refer to \cite{Lieb0,Lions,Lions1,MS,MVs} and references therein for a complete overview). Indeed their solutions are steady states of a generalized nonlinear Schr\"odinger equation, with an attractive interaction potential  given in terms of the Riesz interaction kernel,  which is therefore  weaker and with longer range than  the usual power-type potential in nonlinear Schr\"odinger equation.  We recall that  the relation between MFG systems and normalized nonlinear elliptic equations  has been exploited    in the recent work \cite{PePiVaVe} in the case of nonlinear Schr\"odinger systems.\\

Going back to our Mean-Field Game system, the two distinctive features of our model are the following: the state space is the whole Euclidean space $\R^N$, and the coupling is aggregative and  defined in terms of a Riesz-type interaction kernel. Usually, Mean-Field Game systems are considered in bounded domains, with Neumann or periodic boundary conditions, in order to avoid non-compactness issues. We recall some works in the non compact setting: in particular \cite{BaPi} in the linear-quadratic framework, \cite{Po} in the time-dependent case, \cite{GoPi} for regularity results and finally \cite{CC}, where a system analogous to \eqref{1} has been considered, with power-type nonlinearity. In the unbounded setting, the dissipation induced by the Brownian motion has to be compensated by the optimal velocity, which is a priori unknown and depends by the distribution $m$ itself and on the coercive potential $V$. The  coercive potential $V$ describes spatial preferences of agents and hence discourages them to be far away from the origin. Moreover, due to the presence of the Riesz-type interaction potential $-K_\al\ast m$ which represents the coupling between the individual and the overall population, every player of the game is attracted toward regions where the population is highly distributed. Most of the MFG literature focuses on the study of systems with competition, namely when the coupling descourages aggregation: this  assumption is essential if one seeks for uniqueness of equilibria,
and it is in general crucial in many existence and regularity arguments (see \cite{GoPibook}). Focusing MFG systems, namely models with  coupling which encourages aggregation,  have been studied for instance in \cite{CC,CeCi,Cir,Cir1,GNP} in the stationary setting.\\

In this paper we provide existence and nonexistence results of classical solutions solving the MFG system \eqref{1}, where by \textit{classical solution} we will mean a triple $(u,m,\la)\in C^2(\R^N)\times W^{1,p}(\R^N)\times \R$ for every $p\in(1,+\infty)$.  Our focus will be to obtain classical solutions which satisfy some integrability  conditions and boundary conditions at $\infty$  which will be meaningful from the point of view of the game. 
In particular, we will require some integrability properties of the optimal speed with respect to $m$, namely 
\begin{equation}\label{integr} m|\nabla u|^\g\in L^1(\R^N) \qquad Vm\in L^1(\R^N)\qquad\text{and}\qquad|\nabla m||\nabla u|\in L^1(\R^N).\end{equation}
Indeed, if one looks at the Kolmogorov equation, such integrability properties are important   to ensure some minimal regularity
of $m$  and uniqueness of the invariant distribution itself (see \cite{H, MPR}).  Regularity and boundedness of $m$  is quite crucial in our setting:  indeed, due to the aggregating forces, $m$ has an intrinsic tendency to concentrate and hence to develop singularities. Moreover the Lagrange multiplier $\lambda$ will be uniquely defined as   the generalized principal eigenvalue (see for details \cite{BM, Ci, Ic1}): for fixed  $m\in L^1(\R^N)$ such that $K_\alpha\ast m\in C^{0, \theta}(\R^N)$ for some $\theta\in (0,1)$, 
we define $\la$ as
$$\la:=\sup\left\{c\in\R\,\bigg|\, \exists v\in C^2(\R^N) \text{ solving }\-\Delta v+\frac{1}{\g}|\nabla v|^\g+c=V(x)-K_\al\ast m\right\}.$$
Once that we know this value exists, it is possible to show that there exists $u\in C^2(\R^N)$ solving the HJB equation with such value $\la$, and that such solution $u$ is coercive i.e.  
\begin{equation}\label{crescita} 
u(x)\to +\infty\qquad\text{as $|x|\to +\infty$}
\end{equation}  
and moreover its gradient has polynomial growth (see Section 2 and the references \cite{BM, Ci, Ic1}). Note that \eqref{crescita} is a quite natural “boundary” condition for ergodic HJB equations on the whole space: indeed the optimal speed would give rise to an ergodic process,  so in particular, at least heuristically $-\nabla u \cdot x<0$ for $|x|\to +\infty$, (refer to \cite{H} and references therein, for more information about ergodic problems on the whole space and their characterization in terms of Lyapunov functions). Existence results  for such classical solutions  will depend on the interplay between the dissipation (i.e. by the diffusive term in the system) and the aggregating forces (described in terms of the Riesz potential $K_\al$ and the coercive potential $V$). So, we get that the MFG system \eqref{1} shows three different regimes which correspond to $\al\in(0,N-2\g')$, $\al\in (N- 2\g',N-\g')$ and $\al\in (N-\g',N)$. We will refer to $\al=N-2\g'$ as the \textit{Hardy-Littlewood-Sobolev-critical exponent} and to $\al=N-\g'$  as the \textit{mass-critical} (or $L^2$-critical) \textit{exponent}, in analogy with the regimes appearing in the study of the Choquard equation \eqref{cho} when $\g'=2$. Obviously if $\g'\geq N$, there exists just one regime, which will be the mass-subcritical regime $\alpha \in [0, N)$, whereas if $\frac{N}{2}\leq \g'<N$ there will be just 2 regimes.\\

First of all we observe that for classical solutions to \eqref{1} with $V\equiv 0$ and which satisfy \eqref{integr}, a Pohozaev type identity holds  (see Proposition \ref{prop3.1}): 
\begin{equation}\label{pohointro}
(2-N)\int\limits_{\R^N}\nabla u\cdot\nabla m \,dx+\left(1-\frac{N} {\g}\right)\int \limits_{\R^N} m|\nabla u|^\g dx=\la N M+\frac{\al+N}{2}\int\limits_{\R^{2N}}\frac{m(x)m(y)}{|x-y|^{N-\alpha}}dxdy.
\end{equation}  Also in presence of the potential a similar identity holds, under the  additional integrability condition that $m \nabla V\cdot x \in L^1(\R^N)$. For MFG in the periodic setting with polynomial interaction potential an analogous Pohozaev identity has been proved in \cite{Cir}. For the case of the Choquard equation we refer to \cite{MS} and references therein.   

In the \textit{Hardy-Littlewood-Sobolev-supercritical regime} $0<\al<N-2\g'$, the Pohozaev identity, together with the fact that $\lambda\leq 0$ (see Lemma \ref{lem3.3}),  implies that solutions to the MFG system \eqref{1} do not exist. More precisely, we obtain the following nonexistence result.

\begin{theorem}\label{teo_non_esist}
Assume that $\al\in(0,N-2\g')$ and $V\equiv 0$. Then, the MFG system \eqref{1} has no classical  solutions $(u,m,\la)\in C^2(\R^N)\times W^{1,\frac{2N}{N+\al}}(\R^N)\times\R$ which satisfy \eqref{integr} and \eqref{crescita}.
\end{theorem}

In the case when $N-2\g'<\al<N$ we obtain existence of classical solutions to the MFG system \eqref{1} by means of a Schauder fixed point argument (refer to \cite{BaFe}  and see also \cite{Cir}). More in detail, we consider a regularized  version of problem \eqref{1}, obtained by convolving  the Riesz-interaction term with a sequence of standard symmetric mollifiers (see \eqref{1 phi} below). Taking advantage of the fixed-point structure associated to the MFG system and exploiting  the Schauder Fixed Point Theorem, we show that solutions to the ``regularized'' version of the MFG system do exist. Then we provide a priori uniform estimates on the solutions   to the regularized problem, which allow us to pass to the limit and obtain a classical  solution of the MFG system \eqref{1}.  

\begin{theorem}\label{teo esist fix}
Assume that the potential $V$ is locally H\"older continuous and satisfies \eqref{V}. We have the following results:
\begin{itemize}
\item[i.] if $N-\g'<\al<N$ then, for every $M>0$ the MFG system \eqref{1} admits a classical solution $(u,m,\lambda)$;
\item[ii.] if $N-2\g'<\al\leq N-\g'$ then, there exists a positive real value $M_0=M_0(N,\al,\g,C_V,b)$ such that if $M\in(0, M_0)$ the MFG system \eqref{1} admits a classical solution $(u, m, \lambda)$.
\end{itemize}
Moreover in both cases   there exists a constant $C>0$ such that
$$|\nabla u(x)|\le C(1+|x|)^\frac{b}{\g}\qquad     u(x)\ge C|x|^{\frac{b}{\g}+1}-C^{-1},$$
where $C=C(C_V,b,\g, N,\la,\al)$, $\sqrt{m}\in W^{1,2}(\R^N)$  and it holds
$$m|\nabla u|^\g\in L^1(\R^N),\qquad m V\in L^1(\R^N),\qquad |\nabla u|\,|\nabla m|\in L^1(\R^N).$$
\end{theorem}

Note that in  the mass-subcritical case, solutions to the MFG exist for every mass $M$, whereas in the  mass-supercritical case  and mass-critical case (namely for $\al\in(N-2\g',N-\g']$) we provide existence just for sufficiently small masses, below some threshold $M_0$, due to the fact that in this case the interaction attractive potential is stronger than the diffusive part.  

The Hardy-Littlewood-Sobolev critical exponent is not covered by our analysis. Indeed it is possible to prove existence of solutions to the regularized problem also in this case, for sufficiently small masses (see Theorem \ref{esist MFGreg}). Nevertheless in order to pass to the limit in the regularization, we need to obtain a priori $L^\infty$ bounds on solutions  $m_k$ to the regularized problem, starting from uniform bounds in $L^{\frac{2N}{N+\al}}\cap L^1$. This is not possible at the critical level $\al=N-2\g'$, due to critical rescaling properties of the Sobolev critical exponent:  a priori uniform $L^\infty$ bounds on $m_k$ only holds in the range when we have a uniform bound in $L^q$, for $q>\frac{N}{\g'+\alpha}$ (see Theorem \ref{teo4.1new}) and $\frac{2N}{N+\al}>\frac{N}{\g'+\al}$ only in the  Hardy-Littlewood-Sobolev subcritical  regime. One way to circumvent this difficulty would be to obtain at the critical level $\al=N-2\g'$, by using regularity estimates on the viscous Hamilton-Jacobi equation and on the Fokker Planck equation and a smallness condition on $\|m\|_{\frac{N}{N-\g'}}$, a priori uniform bounds on $m$ in $L^q$ for some $q>\frac{N}{N-\g'}$, in order to be able to apply Theorem \ref{teo4.1new}. This kind of result has been obtained recently in \cite{CirCosVerzini} for MFG in bounded domains with Neumann boundary conditions, and with a nonlinear Schr\"odinger type potential. This problem is related to the maximal regularity of  solutions to viscous Hamilton-Jacobi equation  $-\Delta u+|\nabla u|^\g= f(x)$  (see \cite{cg,cv,g}).  When $m\in L^{\frac{N}{N-\g'}}$, then by Hardy-Littlewood-Sobolev inequality (refer to Theorem \ref{HLS}) $K_\al\ast m\in L^{\frac{N}{\g'}}$, which is a critical threshold in this setting. 
 
To understand better this difference and also the deep analogy with  normalized Choquard-type equations, it will be useful to analyze the problem from a variational point of view.  Existence of solutions to the normalized Choquard equations has been first investigated using variational methods by E.H. Lieb \cite{Lieb} and P.-L. Lions \cite{Lions,Lions84}, while more recently Li-Ye \cite{LiYe}  studied existence of positive solutions to \eqref{cho} by using a minimax procedure and the concentration-compactness principle. As Lasry and Lions first pointed out in \cite{LL}, solutions to \eqref{1} correspond to critical points of the following energy 

\begin{equation}\label{energy}
\mathcal{E}(m,w):=\begin{cases} 
\int\limits_{\R^N}\frac{m}{\g'}\left|\frac{w}{m}\right|^{\g'}+V(x)\,m\,dx-\frac{1}{2}\int\limits_{\R^N}\int\limits_{\R^N} \frac{m(x)\,m(y)}{|x-y|^{N-\al}} dx\,dy\qquad \text{if}\,\,(m,w)\in\mathcal{K}_{M},\\
+\infty \hspace{7.6cm} \text{otherwise}\end{cases}
\end{equation}
where $w:=-m|\nabla u|^{\g-2}\nabla u$ and the constraint set is defined as

\begin{equation}\label{K}
\begin{aligned}
\mathcal{K}_{M}:=\Big\{(m,w)\in & (L^1(\R^N)\cap L^q(\R^N))\times L^1(\R^N)\quad\text{s.t.}\quad \int_{\R^N}m\,dx=M,\quad m\ge0\,\,\, \text{a.e.}\\
&\int_{\R^N}m(-\Delta\p)\,dx=\int_{\R^N}w\cdot\nabla\p\,dx\quad\forall\p\in C^\infty_0(\R^N)\Big\}
\end{aligned}
\end{equation}
with
\begin{equation}\label{q}
q:=\begin{cases}\frac{N}{N-\gamma'+1}\quad\text{if}\,\,\gamma'< N\\
\gamma'\qquad\quad\,\,\,\text{if}\,\,\gamma'\ge N\end{cases}.
\end{equation}
If $N-\g'<\al<N$, so in the \textit{mass-subcritical regime}, the energy $\E$ is bounded from below, indeed using elliptic regularity results for the Kolmogorov equation (see Proposition \ref{prop1.2} below, in this case $1<\frac{2N}{N+\al}<1+ \frac{\g'}{N}$ hence we can use \eqref{stima m_p 1}), the Hardy-Littlewood-Sobolev inequality (see \eqref{hlsenergy}) and the fact that $V\ge0$, we get
$$\E(m,w)\ge C_1 \e^{\g'}\|m\|_{L^\beta(\R^N)}^{\frac{2\g'}{N-\al}}-C_2\|m\|^2_{L^\beta(\R^N)}$$
where we denoted $\beta=\frac{2N}{N+\al}$. Hence, $\inf_{(m,w)\in \mathcal{K}_M} \mathcal{E}(m,w)$ is well defined and by means of classical direct methods and compactness arguments, it is possible to construct global minimizers. Then a linearization argument and a convex duality theorem allow us to show that  minimizers $(m,w)$ of $\E$ correspond to  solutions to the MFG system \eqref{1} (for more details we refer to \cite{CC, CeCi}). In the  \textit{mass-critical regime}, namely for $\al=N-\g'$, the energy is bounded from below just for sufficiently small masses $M$, and we may construct in this range global minimizers. In the \textit{mass-supercritical regime}, namely for $0<\al<N-\g'$, the energy is not bounded from below in general so no global minimum can be found. Nonetheless some compactness of  sequences with finite energy  is still available in the \textit{Hardy-Littlewood-Sobolev-subcritical regime} $N-2\g'\leq \al<N-\g'$. In particular, we may consider a  minimization problem adding a smallness constraint on the $L^\frac{2N}{N+\al}(\R^N)$  norm of $m$ and we may show that if the total mass of $m$ is sufficiently small, then the constrained minimizers are actually local free minimizers of the problem. This procedure would provide solutions to the Mean-Field Game which should coincide with the solutions we  obtained in Theorem \ref{teo esist fix} for $\al\in(N-2\g',N-\g')$  by using Schauder Fixed Point Theorem and imposing a smallness condition on  the $L^\frac{2N}{N+\al}(\R^N)$  norm of $m$.  A similar procedure for constructing local minimizers  has recently been developed for MFG in bounded domains with  Neumann boundary conditions and local aggregative interaction potential of polynomial type (i.e. with a  nonlinear Schr\"odinger type potential), we refer to \cite{CirCosVerzini}.  Moreover, since the energy is becoming more and more negative as the $L^\frac{2N}{N+\al}(\R^N)$ norm of $m$ increases (as it can be observed by a simple rescaling argument), then we expect that with a nontrivial adaptation of the mountain-pass theorem, it should be possible to construct in the Hardy-Littlewood-Sobolev-subcritical regime $N-2\gamma'< \alpha<N-\gamma'$ also solutions to the MFG with a min-max procedure (analogously to what is done in the case of normalized Choquard equation, see \cite{LiYe}). We plan to investigate this issue in a forthcoming paper.

Finally we leave open the problem of existence of classical solutions to the MFG system for $\al\in [N-2\g',N)$ when $V\equiv 0$. Using the   variational approach, and  an appropriate adaptation of the concentration-compactness Lions theorem, one of the author provided existence of solution in the mass subcritical regime $\al\in (N- \g',N)$ as global minimizers of the energy \eqref{energy} with $V\equiv 0$ among competitors with appropriate integrability condition at infinity, see \cite{ber}. We expect that  in the supercritical mass regime  $\al\in (N- 2\g',N-\gamma')$ local minimizers  of the free energy are not present, but that it could be possible to construct a critical point of the energy by means of concentration-compactness arguments together with a min-max procedure. We plan to investigate this issue in a forthcoming paper.
 
The paper is organized as follows. Section 2 contains some preliminary results. In particular we recall  regularizing properties of the Riesz interaction kernel, some a priori elliptic  estimates for solutions to the Kolmogorov equation,  a priori gradient estimates for solutions to the Hamilton Jacobi Bellman equation and finally   uniform $L^\infty$  bounds for $m$, solution to \eqref{1}. 
In Section 3 we provide the Pohozaev identity and the proof of the non-existence result Theorem \ref{teo_non_esist}. Section 3 contains the proof of the existence result Theorem \ref{teo esist fix}. 

In what follows, $C, C_1, C_2, K_1,\dots$ denote generic positive constants which may change from line to line and also within the same line. 
Moreover $\gamma'$ denotes the conjugates exponent of $\gamma$, that is $\gamma'=\frac{\gamma}{\gamma-1}$.


\section{Preliminaries}\label{sez2}

\subsection{Regularity results for the Kolmogorov equation}

\begin{lemma}\label{mC2}
Let $u\in C^{2,\theta}(\R^N)$ and $m\in W^{1,2}(\R^N)$ be a solution (in the distributional sense) to 
\begin{equation}\label{10}
-\Delta m(x)-\mathrm{div}\left(m(x)\,\nabla u(x)\,|\nabla u|^{\gamma-2}\right)=0 \quad \text{in}\,\,\R^N,
\end{equation}
where $\g>1$ is fixed. Then, $m\in C^{2,\alpha}(\R^N)$. Moreover, if $m\ge0$ and  $m\not\equiv 0$,  then $m(x)>0$ for any $x\in\R^N$.
\end{lemma}

\proof 
If $\g\ge2$, then $m$ solves 
$$-\Delta m-b(x)\cdot\nabla m(x)-m(x)\,\mathrm{div}\,b(x)=0$$
where $b(x):=|\nabla u|^{\g-2}\nabla u(x)\in C^{1,\theta}(\R^N)$ and $\mathrm{div}\,b(x) \in C^{0,\theta}(\R^N)$. By elliptic regularity (see e.g. \cite[Theorem 8.24]{GT}) we get that $m\in C^{0,\al}$. Denoting by $f:=m\nabla u|\nabla u|^{\g-2}$ we have $-\Delta m=\mathrm{div} f$ where $f\in C^{0,\al}$, then by \cite[Theorem 4.15]{GT} we get that $m\in C^{1,\al}$ and hence 
$$-\Delta m=\mathrm{div}\left(m\nabla u |\nabla u|^{\g-2}\right)\in C^{0,\min\{\al,\theta\}}$$
so $m\in C^{2,\min\{\al,\theta\}}$, and  iterating we finally obtain  that $m\in C^{2,\theta}$. If $1<\g<2$, $b(x)$ is just an H\"older continuous function, hence $m$ is a weak solution of equation \eqref{10}. In this case, we can replace $b(x)$ with $b_\varepsilon(x):=\nabla u(x)(\varepsilon+|\nabla u|^2)^{\frac{\gamma}{2}-1}$ and $m_\varepsilon$ is a-posteriori a classical solution to the approximate equation $$-\Delta m-\mathrm{div}\left(m(x)\,b_\varepsilon(x)\right)=0.$$
We can conclude letting $\varepsilon\to0$. If $m\ge0$ on $\R^N$, we also have that $m$ satisfy
$$-\Delta m-b(x)\cdot\nabla m(x)-\big(\mathrm{div}\,b(x)\big)^+m(x)\le0,$$
since $\int_{\R^N}m\,dx=M>0$, the Strong Minimum Principle (refer e.g to \cite[Theorem 8.19]{GT}) implies that $m>0$ in $\R^N$ (indeed $m$ can not be equal to 0, unless it is constant, which is impossible).
\endproof

We will use the following result (proved in \cite[Proposition 2.4]{CC}) which takes advantage of some classical elliptic regularity results of Agmon \cite{Ag}.

\begin{prop}\label{prop2.4CC}
Let $m\in L^p(\R^N)$ for $p>1$ and assume that for some $K>0$ 
$$\left|\int_{\R^N}m\Delta\varphi\,dx\right|\le K\|\nabla\varphi\|_{L^{p'} (\R^N)}, \quad \forall\varphi\in C_0^\infty(\R^N).$$
Then, $m\in W^{1,p}(\R^N)$ and there exists a constant $C>0$ depending only on $p$ such that 
$$\|\nabla m\|_{L^p(\R^N)}\le C\,K.$$
\end{prop}

We prove now some a priori estimates for solutions to the Kolmogorov equation. Let us fix $p\in(1,+\infty)$ and $M>0$.

\begin{prop}\label{prop1}
Let us consider a couple $(m,w)\in (L^p(\R^N)\cap W^{1,1}(\R^N))\times L^1(\R^N)$ which solves weakly
$$-\Delta m+\mathrm{div}\,w=0,\quad in\,\,\,\R^N.$$
Assume also that $\int_{\R^N}m(x)\,dx=M$, $m\ge0$ a.e. and 
\begin{equation}\label{c}
E:=\int_{\R^N}m\left|\frac{w}{m}\right|^{\g'}dx<+\infty.
\end{equation}
Then, we have that
$$m\in W^{1,r}(\R^N)$$
for $r$ such that $\frac{1}{r}=\left(1-\frac{1}{\g'}\right)\frac{1}{p}+\frac{1}{\g'}$ (i.e. $r=\frac{p\g'}{\g'+p-1}$) and there exists a constant $C$, depending on $r$, such that
\begin{equation}\label{mW}
\|m\|_{W^{1,r}(\R^N)}\le C(E+M)^{\frac{1}{\gamma'}}\|m\|_{L^p(\R^N)} ^{\frac{1}{\gamma}}.
\end{equation}
\end{prop}

\proof
By definition of weak solution we have
$$-\int_{\R^N}m\,\Delta\varphi\, dx=\int_{\R^N}w\cdot\nabla \varphi\,dx,\quad\text{for every} \,\,\,\varphi\in C_0^\infty(\R^N),$$
using H\"older inequality (since $\frac{1}{r}=\left(1-\frac{1}{\g'}\right)\frac{1}{p}+ \frac{1}{\g'}$, it holds $\frac{1}{p\g}+\frac{1}{\g'} +\frac{1}{r'}=1$) we obtain
\begin{align*}
\left|\int_{\R^N}m\,\Delta\varphi\, dx\right|&\le\int_{\R^N}|w|\,|\nabla \varphi|dx=\int_{\R^N}\left(\left|\frac{w}{m}\right|^{\g'} m\right)^ {\frac{1}{\g'}} m^{\frac{1}{\g}} |\nabla\varphi|\,dx\\ &\le\left(\int_{\R^N}\left|\frac{w}{m}\right|^{\g'}m\,dx\right)^ {\frac{1}{\g'}}\|m\|_{L^p(\R^N)} ^{\frac{1}{\g}}\|\nabla\varphi\|_{L^{r'}(\R^N)}
\end{align*}
and hence
$$\left|\int_{\R^N}m\,\Delta\varphi\, dx\right|\le E^ {\frac{1} {\g'}}\|m\|_{L^p(\R^N)}^{\frac{1}{\g}}\|\nabla\varphi\|_{L^{r'}(\R^N)}.$$
Since $\|m\|_{L^1(\R^N)}=M$ and $m\in L^p(\R^N)$, by interpolation we get 
\begin{equation}\label{mLr}
\|m\|_{L^r(\R^N)}\le\|m\|_{L^p(\R^N)}^{\frac{1}{\g}}M^{\frac{1}{\g'}}
\end{equation}
therefore $m\in L^r(\R^N)$. From Proposition \ref{prop2.4CC} with $K=E^{\frac{1}{\gamma'}}\|m\| _{L^p(\R^N)}^{\frac{1}{\gamma}}$, we obtain that $m\in W^{1,r}(\R^N)$ and there exists a constant $C>0$, depending on $r$, such that 
\begin{equation}\label{DmLr}
\|\nabla m\|_{L^r(\R^N)}\le C\,E^{\frac{1}{\gamma'}}\|m\| _{L^p(\R^N)}^{\frac{1}{\gamma}}.
\end{equation}
By \eqref{mLr} and \eqref{DmLr}, we can conclude that
$$\|m\|_{W^{1,r}(\R^N)}\le \left(M^{\frac{1}{\gamma'}}+CE^{\frac{1}{\gamma'}} \right)\|m\|_{L^p}^ {\frac{1}{\gamma}}\le C(E+M)^{\frac{1}{\gamma'}}\|m\|_{L^p(\R^N)} ^{\frac{1}{\gamma}}.$$
\endproof

\begin{prop}\label{prop1.2}
Under the assumption  of Proposition \ref{prop1}, we have the following results:
\begin{itemize}
\item[i)] if $1<p<1+\frac{\g'}{N}$ then, there exists $\delta_1=\frac{1}{p-1}\left (\frac{\g'}{N}+1-p\right)$ such that
\begin{equation}\label{stima int w/m}
    \|m\|_{L^p(\R^N)}^{(1+\delta_1)p}\le C\,M^{(1+\delta_1)p-1}\,E
\end{equation}
where $C$ is a constant depending on $N$, $\g$ and $p$;
\item[ii)] if $\g'<N$ and $1<p\le\frac{N}{N-\g'}$ then, there exists $\delta_2=\frac{1}{p-1}\frac{\g'}{N}$ and a constant $C$ depending on $N$, $\g$ and $p$ such that
\begin{equation}\label{mdeltap}
    \|m\|_{L^p(\R^N)}^{p\delta_2}\le C(E+M)M^{p\delta_2-1}.
\end{equation}
\end{itemize}
\end{prop}

\proof
$i)$ The proof of \eqref{stima int w/m} follows from \cite[Lemma 2.8]{CC}. $ii)$ As before let $\frac{1}{r}=\frac{1}{p}\left(1-\frac{1}{\g'}\right)+\frac{1}{\g'}$, if $\g'<N$ then $r<\g'<N$, so by Gagliardo-Niremberg inequality and \eqref{mW} we get
\begin{equation}\label{d}
\|m\|_{L^{r^*}(\R^N)}\le C\| m\|_{L^p(\R^N)}^{\frac{1}{\g}} (E+M)^{\frac{1}{\g'}}  
\end{equation}
where $\frac{1}{r^*}=\frac{1}{r}-\frac{1}{N}$ and $C$ is a constant depending on $N$, $p$ and $\g'$. One can observe that $\frac{1}{r^*}-\frac{1}{p}= \frac{p N-N-p\g'}{p\g' N}\le0$, that is $r^*\ge p$, by interpolation there exists $\theta\in(0,1]$ such that 
$$\|m\|_{L^p(\R^N)}^{\frac{1}{\theta}}\le M^{\frac{1-\theta}{\theta}} \|m\|_{L^{r^*}(\R^N)}$$
and from \eqref{d} we get that
$$\|m\|_{L^p(\R^N)}^{\left(\frac{1}{\theta}-\frac{1}{\g}\right)\g'}\le C (E+M) M^{\frac{1-\theta}{\theta}\g'}.$$
By simple computations we have that 
$$\left(\frac{1}{\theta}-\frac{1}{\g}\right)\g'=\frac{\g'}{N}\frac{p}{p-1}$$
and 
$$\left(\frac{1}{\theta}-1\right)\g'=\frac{\g'}{N}\frac{p}{p-1}-1$$
denoting by $\delta_2$ the quantity $\frac{1}{p-1}\frac{\g'}{N}$, we finally obtain \eqref{mdeltap}.
\endproof

\begin{obs}
In the following we will use \eqref{stima int w/m} and \eqref{mdeltap} in the case when $p=\frac{2N}{N+\al}$. It will be useful to observe that if $\g'\ge N$ then $1<\frac{2N}{N+\alpha}<2\le 1+\frac{\g'}{N}$, hence estimate \eqref{stima int w/m} holds. In the case when $\g'<N$, if $N-\g'\le\al<N$ then, $1<\frac{2N}{N+\alpha}<1+\frac{\g'}{N}$ and hence from \eqref{stima int w/m} we get that
\begin{equation}\label{stima m_p 1}
\|m\|_{L^\frac{2N}{N+\al}(\R^N)}^\frac{2\g'}{N-\al}\le CM^{\frac{2\g'}{N-\al}-1}E;
\end{equation}
whereas if $N-2\g'\le\al<N-\g'$, we may use estimate \eqref{mdeltap}, which gives us
\begin{equation}\label{stima m_p 2}
\|m\|_{L^\frac{2N}{N+\al}(\R^N)}^\frac{2\g'}{N-\al}\le C(E+M) M ^{\frac{2\g'}{N-\al}-1}.
\end{equation}
\end{obs}

Finally, we recall the following a priori elliptic regularity result (see \cite[Proposition 2.8, Corollary 2.9]{CC}). 

\begin{prop}\label{kolm}
Let \[ q:=\begin{cases}\frac{N}{N-\g'+1}\quad\text{if}\,\,\g'< N\\
\g'\qquad\quad\,\,\,\text{if}\,\,\g'\ge N\end{cases}
\]  and let   $(m,w)\in (L^q(\R^N)\cap L^1(\R^N))\times L^1(\R^N)$ be  a weak solution to 
$$-\Delta m+\mathrm{div}\,w=0,\quad in\,\,\,\R^N$$
with  $\int_{\R^N}m(x)\,dx=M$, $m\ge0$ a.e. and 
\[
E:=\int_{\R^N}m\left|\frac{w}{m}\right|^{\g'}dx<+\infty.
\]
Then the following hold:
\begin{enumerate}
    \item[i)] $$m\in L^\beta(\R^N),\quad\forall \beta\in\left[1,\frac{N}{N-\g'}\right)\qquad(\forall\beta\in[1,+\infty), \,\,\,\text{if}\,\,\, \g'\ge N)$$
    and there exists a constant $C$ depending on $N$, $\beta$ and $\gamma'$ such that $$\|m\|_{L^\beta(\R^N)}\le C(E+M);$$
    \item[ii)] $$m\in W^{1,\ell}(\R^N),\quad\forall \ell<q$$
    and there exists a constant $C$ depending on $N$, $\ell$ and $\g'$ such that $$\|m\|_{W^{1,\ell}(\R^N)}\le C(E+M).$$
\end{enumerate}
\end{prop}

\proof
From Proposition \ref{prop1} we have
$$m\in W^{1,r_0}(\R^N)\quad\text{for}\quad \frac{1}{r_0}=\left(1-\frac{1}{\g'}\right) \frac{1}{q}+\frac{1}{\g'}.$$
\indent \textit{Case $\g'<N$}. Since $1<r_0<\g'<N$, by Sobolev embedding theorem and interpolation, we get that 
\begin{equation}\label{sobemb}
m\in L^{\beta}(\R^N)\quad \forall \beta\le q_1
\end{equation}
where $q_1$ is the Sobolev critical exponent, i.e.
$$q_1:=\frac{Nr_0}{N-r_0}=\frac{qN\g'}{N\g'-N+q(N-\g')},$$
(notice that $q_1>q$ since $q<\frac{N}{N-\gamma'}$). From \eqref{sobemb}, using Proposition \ref{prop1} again, we have
$$m\in W^{1,\ell}(\R^N)\quad\quad \forall \ell\le r_1=\frac{q_1 \g'}{\g'-1+q_1}.$$
As before, by Sobolev embedding theorem and interpolation, we have that 
$$m\in L^{\beta}(\R^N)\qquad\forall\beta\le q_2=\frac{q_1 N\g'}{N\g'-N+q_1(N-\g')}.$$
Iterating the previous argument, we observe that $q_{j+1}=f(q_{j})$ where $f(s):=\frac{sN\g'}{N\g'-N+s(N-\g')}$. Since $f$ is an increasing function if $s<\frac{N}{N-\g'}$ and it has a fixed point for $\Bar{s}=\frac{N}{N-\g'}$, we obtain that
$$m\in L^\beta(\R^N),\quad \forall\beta<\frac{N}{N-\g'}$$
and
$$m\in W^{1,\ell}(\R^N),\quad \forall\ell<\frac{N}{N-\g'+1}.$$
Moreover, for any fixed $\beta<\frac{N}{N-\g'}$, taking $r=r(\beta)$ such that $\frac{1}{r}=\left(1-\frac{1}{\g'}\right)\frac{1}{\beta}+\frac{1}{\g'}$, from estimate \eqref{mW} and the Sobolev embedding theorem (notice that $r^*>\beta$) we
get that there exists a constant $C$ depending on $N$ and $r$ such that
\begin{equation}\label{norm}
\|m\|_{L^\beta(\R^N)}\le C(E+M)^{\frac{1}{\g'}}\|m\|_{L^\beta(\R^N)}^ \frac{1}{\g}.
\end{equation}
and hence
\begin{equation}\label{norm m}
\|m\|_{L^\beta(\R^N)}\le C_1(E+M).
\end{equation}
Putting \eqref{norm m} in \eqref{mW} we obtain
$$\|m\|_{W^{1,\ell}(\R^N)}\le C_2(E+M).$$
\indent \textit{Case $\g'=N$.}
Since $r_0<\g'=N$, we can apply the Sobolev embedding theorem and with the same argument as before we obtain
$$q_{j+1}=\frac{N}{N-1}q_j.$$
Obviously $q_{j+1}>q_j$, by iteration we get that 
$$m\in L^\beta(\R^N),\quad \forall\beta<+\infty$$
and 
$$m\in W^{1,\ell}(\R^N),\quad \forall\ell<\g'.$$
The estimates on the norms follow in the same way as the previous case.\\
\indent \textit{Case $\g'>N$.} Since $m\in L^{\g'}(\R^N)$, by interpolation $m\in L^N(\R^N)$ and we can go back to the previous case.
\endproof





\subsection{Some properties of the Riesz potential}

We recall here some properties of the Riesz potential, which will be useful in the following in order to deal with the Riesz-type interaction term. 

\begin{definition}
Given $\al\in(0,N)$ and a function $f\in L^1_{loc}(\R^N)$, we define the Riesz potential of order $\al$ of $f$ as
$$K_\al\ast f (x):=\int_{\R^N}\frac{f(y)}{|x-y|^{N-\al}}dy,\quad x\in\R^N.$$
\end{definition}

The Riesz potential $K_\alpha$ is well-defined as an operator on the whole space $L^r(\R^N)$ if and only if $r\in\left[1,\frac{N}{\alpha}\right)$. We state now the following well-known theorems (for which refer e.g. to \cite[Theorem 14.37]{WZ} and \cite[Theorem 4.3]{LL}).

\begin{theorem}[Hardy-Littlewood-Sobolev inequality] \label{HLS}
Let $0<\al<N$ and $1<r<\frac{N}{\al}$. Then for any $f\in L^r(\R^N)$ 
$$\|K_\alpha\ast f\|_{L^{\frac{Nr}{N-\alpha r}}(\R^N)}\le C\|f\|_{L^r(\R^N)}$$
where $C$ is a constant depending only on $N$, $\alpha$ and $r$.
\end{theorem}

\begin{theorem}\label{HLS2}
Let $0<\lambda<N$ and $p,r>1$ with $\frac{1}{p}+\frac{\lambda}{N} +\frac{1}{r} =2$. Let $f\in L^p(\R^N)$ and $g\in L^r(\R^N)$. Then, there exists a sharp constant $C(N,\lambda,p)$ (independent of $f$ and $g$) such that
\begin{equation}\label{disug hls}
\left|\int_{\R^N}\int_{\R^N}\frac{f(x)\,g(y)}{|x-y|^{\lambda}}dx\,dy\right|\le C\|f\|_{L^p(\R^N)}\|g\|_{L^r(\R^N)}.
\end{equation}
\end{theorem}

\begin{obs}
It follows immediately that if $0<\al<N$ and $f\in L^\frac{2N}{N+\alpha}(\R^N)$, then there exists a sharp constant $C$, depending only on $N$ and $\alpha$, such that
\begin{equation}\label{hlsenergy}
\left|\int_{\R^N}\int_{\R^N}\frac{f(x)\,f(y)}{|x-y|^{N-\alpha}}dx\,dy\right|\le C\|f\|^2_{L^\frac{2N}{N+\alpha}(\R^N)}.
\end{equation}
As shown in \cite{Lieb}, in this case the constant $C$ can be computed explicitly and there exist explicit optimizers for \eqref{hlsenergy} (while neither the constant nor the optimizers are known for $p\not=r$, although do exist).
\end{obs}

Regarding the $L^\infty$-norm and the H\"older continuity of the Riesz potential, we recall  here the following results.

\begin{theorem}\label{holderRiesz}
Let $0<\al<N$,  $1<r\le+\infty$ be such that $r>\frac{N}{\al}$ and $s\in\left[1,\frac{N}{\al}\right)$. Then, for every $f\in L^s(\R^N)\cap L^r(\R^N)$ we have that 
\begin{equation}\label{K infty}
   \|K_\al\ast f\|_{L^\infty(\R^N)}\le C_1\|f\|_{L^r(\R^N)}+C_2\|f\|_{L^s(\R^N)} 
\end{equation}
where $C_1=C_1(N,\al,r)$ and $C_2=C_2(N,\al,s)$. Moreover, if $0<\al-\frac{N}{r}<1$ then, 
$$K_\al\ast f\in C^{0,\al-\frac{N}{r}}(\R^N)$$
and there exists a constant $C$, depending on $r,\,\al$ and $N$, such that  
$$\frac{\big|K_\al\ast f(x)-K_\al\ast f(y)\big|}{|x-y|^{\al-\frac{N}{r}}}\le C \|f\|_{L^r(\R^N)}\qquad \text{ for $x\neq y$}.$$
\end{theorem}

\proof 
Concerning H\"older regularity results for the Riesz potential, one may refer to \cite[Theorem 2.2, p.155]{Mi} and \cite[Theorem 2]{Du}. \\
\indent As for \eqref{K infty}, notice that
$$\frac{1}{|x|^{N-\al}}\in L^p(B_1),\quad \forall p\in\left[1,\frac{N}{N-\al}\right)$$ where $B_1$ is the ball of radius $1$ centered at $0$ 
and it is well-known that $\int_{B_1(0)}\frac{1}{|x|^{(N-\al)p}}dx=\frac{\omega_N}{N-(N-\al)p}$.
By Holder inequality we get
\begin{align*}
\int_{B_1}\frac{|f(x-y)|}{|y|^{N-\al}}dy&\le\bigg(\int_{B_1}|f(x-y)|^r\,dy\bigg)^{\frac{1}{r}}\bigg(\int_{B_1}\frac{1}{|y|^{(N-\al)r'}}\,dy\bigg)^{\frac{1}{r'}} \\
&\le \|f\|_{L^r(\R^N)}\bigg(\frac{\omega_N}{N-(N-\al)r'}\bigg)^\frac{1}{r'}\le C_1\|f\|_{L^r(\R^N)}
\end{align*}
using the fact that $r'<\frac{N}{N-\al}$, since by assumption $r>\frac{N}{\al}$. On the other hand 
$$\frac{1}{|x|^{N-\al}}\in L^p(B_1^c),\quad \forall p\in\left(\frac{N}{N-\al},+\infty\right]$$
hence
$$\int_{\R^N\setminus B_1}\frac{|f(x-y)|}{|y|^{N-\al}}dy\le\left(\int_{\R^N\setminus B_1}|f(x-y)|^s dy\right)^{\frac{1}{s}}\bigg\|\frac{1}{|y|^{N-\al}}dy\bigg\|_{L^{s'}(B_1^c)}\le C\|f\|_{L^{s}(\R^N)},$$
since $(N-\al)s'>N$.
\endproof







\subsection{Some results on the Hamilton-Jacobi-Bellman equation}

By a straightforward adaptation of \cite[Theorem 2.5 and Theorem 2.6]{CC}, we obtain some a priori regularity estimates for solutions to some Hamilton-Jacobi-Bellman equations defined on the whole euclidean space $\R^N$. The following propositions are stated under slightly more general assumptions than ones of our problem.

\begin{prop}\label{lemma2.5}
Assume that $K_\al\ast m\in L^\infty(\R^N)$ and  that $V$ satisfies \eqref{V}, with $b\geq 0$. Let $(u,c)\in C^2(\R^N)\times\R$ be a classical solution to the HJB equation
\begin{equation}\label{HJB}
-\Delta u+\frac{1}{\g}|\nabla u(x)|^\g+c=V(x)-K_\al\ast m(x) \quad in\,\,\R^N,
\end{equation}
for $\g>1$ fixed. Then
\begin{itemize}
\item [i.] there exists a constant $C_1>0$, depending on $C_V,b,\g,N,c,\|K_\al\ast m\|_\infty$, such that 
$$|\nabla u(x)|\le C_1(1+|x|)^{\frac{b}{\gamma}};$$
\item [ii.] if $u$ is bounded from below and $b\not=0$ in \eqref{V}, then there exist a constant $C_2>0$ such that 
$$u(x)\ge C_2|x|^{\frac{b} {\gamma}+1}-C_2^{-1},\quad\forall x\in\R^N.$$
The same result holds also in the case when $b=0$, but we have to require in addition that there exists $\delta>0$ such that $V(x)-K_\al\ast m(x)-c>\delta>0$ for $|x|$ sufficiently large.
\end{itemize}
\end{prop}

\proof
The thesis follows applying \cite[Theorem 2.5 and Theorem 2.6]{CC}.
\endproof

Let us define 
\begin{equation}\label{lambda bar}
\lambda:=\sup\{c \in\R\,|\,\eqref{HJB}\,\,\,\text{has a solution}\,\,\,u\in C^2(\R^N)\}
\end{equation}

\begin{prop}\label{teo2.7}
Besides  the hypothesis of Proposition \ref{lemma2.5}, let us assume also that $V-K_\al\ast m$ is locally H\"older continuous. Then
\begin{itemize}
    \item[i)] $\lambda<+\infty$ and there exists $u\in C^2(\R^N)$ such that the pair $(u, \lambda)$ solves \eqref{HJB}. 
    \item[ii)] if $b\not=0$ in \eqref{V}, $u$ is unique up to additive constants (namely if $(v,\la)\in C^2(\R^N)\times\R$ solves \eqref{HJB} then there exists $k\in\R$ such that $u=v+k$) and there exists a constant $K>0$ such that
    $$u(x)\ge K|x|^{\frac{b} {\g}+1}-K^{-1}, \quad\forall x\in\R^N.$$
\end{itemize}
\end{prop}

\proof 
It follows by \cite[Theorem 2.7]{CC}. We may observe also that \[\lambda=\sup\{c\in\R\, |\,\eqref{HJB}\,\,\,\text{has a subsolution}\,\,\,u\in C^2(\R^N)\}.\]
\endproof 

Finally, we conclude with an estimate on the Lagrange multiplier $\la$ defined in \eqref{lambda bar}. 

\begin{lemma}\label{lem3.3} 
Let $(u,\lambda)\in C^2(\R^N) \times\R$ be a solution to the HJB equation \eqref{HJB}. 
Then \begin{itemize}
\item[i] if $V\equiv 0$, then $\lambda \leq 0$;
\item[ii] if $V$  satisfies \eqref{V} then $\la\leq C$ for some constant depending on $b, C_V,\g,N$. 
\end{itemize} 
\end{lemma}

\proof
The proof is based on the same argument of \cite[Lemma 3.3]{Cir}.
Let us consider the function  $\mu_\delta(x):=\left(\frac{\delta} {2\pi}\right)^{N/2}e^{\frac{-\delta|x|^2}{2}}$ for $x\in\mathbb{R}^N$ and $\delta>0$. Obviously $\int_{\mathbb{R}^N}\mu _\delta(x)\,dx=1$. From the definition of Legendre transform we get that
$$\frac{1}{\g}|\nabla u|^\g=\sup\limits_{\al\in\R^N}\left(\nabla u\cdot \al-\frac {|\al|^{\g'}}{\g'}\right)\ge \nabla u\cdot(\delta x)-\frac {|\delta x|^{\g'}}{\g'}$$
hence
\begin{equation}\label{eq legendre}
-\Delta u(x)+\nabla u\cdot (\delta x)-\frac{1}{\g'}|\delta x|^{\g'}+\la\le V(x)-m\ast K_\al(x).
\end{equation}
Multiplying \eqref{eq legendre} by $\mu_\delta$ and integrating over $B_R$ we obtain
$$-\int_{B_R}\Delta u(x)\mu_\delta+\int_{B_R}\nabla u\cdot (\delta x)\mu_\delta
-\int_{B_R}\frac{1}{\g'}|\delta x|^{\g'}\mu_\delta+\la\int_{B_R}\mu_\delta \le\int_{B_R} (V(x) -m\ast K_\al)\mu_\delta.$$
Integrating by parts (notice that $\int_{B_R}\nabla u\cdot\nabla\mu_\delta=-\int_{B_R} \nabla u\cdot(\delta x)\mu_\delta$) we get
$$-\int_{\partial B_R}\mu_\delta \nabla u \cdot \nu\,d\sigma-\frac{1}{\g'}\int_{B_R}|\delta x|^{\g'}\mu_\delta\,dx+\la\int_{B_R} \mu_\delta\,dx\le\int_{B_R}(V(x)-m\ast K_\al)\mu_\delta\,dx$$
and since $\int_{B_R}m\ast K_\al(x)\mu_\delta(x)\,dx\ge0$, we have 
\begin{equation}\label{45}
\la\int_{B_R} \mu_\delta\,dx\le \int_{\partial B_R}\mu_\delta \nabla u \cdot \nu\,d\sigma+\frac{1}{\g'}\int_{B_R}|\delta x|^{\g'}\mu_\delta\,dx+ \int_{B_R} V(x)\mu_\delta\,dx.
\end{equation}
For $\delta>0$ fixed, the first integral in the RHS of \eqref{45} can be estimated as follows
$$\bigg|\int_{\partial B_R}\mu_\delta \nabla u \cdot \nu\,d\sigma\bigg|\le  C\delta^\frac{N}{2}e^{-\frac{\delta R^2}{2}}\|\nabla u\|_{L^\infty(\partial B_R)}|\partial B_R|\to0,\qquad \text{as}\quad R\to+\infty$$
using  the gradient estimates on $\nabla u$ proved in Proposition \ref{lemma2.5}. 
So, sending $R\to +\infty$ in \eqref{45}  and using \eqref{V} we get 
\[\lambda\leq\frac{1}{(2\pi)^\frac{N}{2}}\frac{\delta^\frac{\g'}{2}}{\g'}\int_{\R^N} |y|^{\g'}e^{-\frac{|y|^2}{2}}\,dy +\frac{1}{(2\pi)^\frac{N}{2}} \int_{\R^N}  V\left(\frac{y}{\sqrt{\delta}}\right) e^{-\frac{|y|^2}{2}}\,dy. \]
If $V\equiv 0$, then sending $\delta\to 0$ in the previous inequality, we conclude immediately $\lambda\leq 0$. 
If $V\not \equiv 0$, we may choose $\delta=1$ in the previous inequality and conclude recalling \eqref{V}. 
\endproof


\subsection{Uniform a priori $L^\infty$-bounds on $m$}

We state now the following result, which provides uniform a priori $L^\infty$ bounds on $m$. 

\begin{theorem}\label{teo4.1new}
  We consider a sequence of classical solutions $(u_k,m_k,\la_k)$ to the following MFG system
\begin{equation}\label{MFGk}
\begin{cases}
-\Delta u+\frac{1}{\g}|\nabla u|^{\g}+\lambda= W_k(x)-G_{k,\al}[m](x)\\
-\Delta m-\mathrm{div} \left(m\nabla u |\nabla u|^{\g-2} \right)=0\\
\int_{\R^N}m=M,\quad m\ge0
\end{cases}\text{in}\,\,\,\R^N
\end{equation}
where $W_k:\R^N\to \R$ satisfies assumption \eqref{V} with constant $C_V,b$ independent of $k$. Let $G_{k,\al}:L^1(\R^N)\to L^1(\R^N)$ such that  $G_{k,\alpha}[m]\geq 0$ for all $m\in L^1$, with $m\geq 0$, and moreover  we assume that there exists $\alpha\in (0, N)$ such that  for all  $s\in\left[1,\frac{N}{\al}\right)$ and $r\in\left(\frac{N}{\al},+\infty\right]$ there holds for $m\in L^s(\R^N)\cap L^r(\R^N)$ 
\begin{equation}\label{stimagk}\|G_{k,\al}[m]\|_{L^\infty(\R^N)}\le C_1\|m\|_{L^r(\R^N)}+C_2\|m\|_{L^s(\R^N)}\end{equation}
where $C_1=C_1(N,\al,r)$ and $C_2=C_2(N,\al,s)$. 

If  $u_k$ are bounded  from below and satisfy \eqref{crescita},  and $m_k\in L^1(\R^N)\cap L^\infty(\R^N)$, with $\|m_k\|_{L^q}\le C_q$ for some $q>\frac{N}{\alpha+\gamma'}$ then, there exists a positive constant $C$ not depending on $k$ such that 
$$\|m_k\|_{L^\infty(\R^N)}\le C, \quad \forall k\in\mathbb{N}.$$
\end{theorem}

\proof
We follow the argument of the proof of \cite[Theorem 4.1]{CC} (refer also to \cite{Cir} for the analogous result on $\mathbb{T}^N$), but we have to define a different rescaling in this case. 

Up to addition of constants we may assume $\inf u_k(x)=0$. 

We assume by contradiction that 
$$\sup_{\R^N} m_k=L_k\to+\infty$$
and we define 
$$\de_k:=\begin{cases} 
L_k^{-\beta}, \,\,\, \quad\text{if }\,\,\g'\le N\quad \text{and } q\leq \frac{N}{\gamma'} \\
L_k^{-\frac{1}{\g'}}, \quad\text{if either }\,\,\g'>N \quad \text{ or } \gamma'\leq N, q>\frac{N}{\gamma'}
\end{cases}
$$
where $\beta>0$ (so $\delta_k\to 0$)  has to be chosen in the following way. We fix \begin{equation}\label{erre} r\in \left(\frac{N}{\alpha}, \frac{Nq}{N-q\gamma'}\right).\end{equation} 
Note that since $q>\frac{N}{\gamma'+\alpha}$ the interval is not empty.  If $q=\frac{N}{\gamma'}$ it is sufficient to choose $\frac{N}{\alpha}<r<+\infty$.  Then we choose $\beta$ such that 
\[\frac{1}{\gamma'}\left(1-\frac{q}{r}\right)\leq \beta< \frac{q}{N}. \]

We rescale $(u_k, m_k, \lambda_k)$ as follows: 
\[v_k(x):=\de_k^{\frac{2-\g}{\g-1}}u_k(\de_k x)+1,\qquad n_k(x):=L_k^{-1}m_k(\de_k x),\qquad \tilde{\la}_k:=\de_k^{\g'}\la_k.\]
Observe that $0\le n_k(x)\le1$ and $\sup n_k=1$ and moreover that $v_k(x)\geq 1$ for all $x$.
So we obtain that $(v_k, n_k, \tilde\la_k)$ is a solution to 
$$\begin{cases}
-\Delta v_k+\frac{1}{\g}|\nabla v_k|^\g+\tilde{\la}_k=V_k(x)-\tilde{g}_k(x)\\
-\Delta n_k-\mathrm{div}(n_k\nabla v_k|\nabla v_k|^{\g-2})=0
\end{cases}$$
where 
\[V_k(x):=\de_k^{\g'} W_k(\de_k x)\qquad\text{and}\qquad \tilde g_k(x):= \de_k^{\g'} G_{k,\al}[m_k](\de_k x).\]
Observe that by assumption \eqref{V} there holds 
\[C_V^{-1}\delta_k^{\gamma'}(\max\{|\delta_k x|-C_V,0\})^b\le V_k(x)\le C_V(1+\delta_k^{\gamma'+b}|x|)^b,\quad \forall x\in\R^N.\]
Computing the equation in a minimum point of $u_k$ we obtain 
$\lambda_k\geq -\|G_{k,\alpha}[m_k]\|_\infty$
and reasoning as in Lemma \ref{lem3.3}, we get that $\lambda_k\leq C$, for some $C$ just depending on $\gamma, C_V, b$,
so we get 
\[-\|\tilde g_k\|_\infty=-\delta_k^{\gamma'} \|G_{k,\alpha}[m_k]\|_\infty\leq \tilde \lambda_k\leq \delta_k^{\gamma'} C. \]

\noindent {\bf If $\gamma'>N$ or $\gamma'\leq N$ and $q> \frac{N}{\gamma'}$ } we apply \eqref{stimagk} with $r=+\infty$ and $s=1$ and we get
\[\|\tilde g_k\|_\infty\leq \de_k^{\g'}(C_1 L_k+C_2M)=L_k^{-1}(C_1 L_k+C_2M)\leq C\] which in turns gives also that $|\tilde \lambda_k|\leq C$. 
Moreover if $\gamma'>N$ there holds  \[  \|n_k\|_{L^1}=\int_{\R^N}n_k(x)dx=\de_k^{\g'-N}\|m_k\|_{L^1}=\de_k^{\g'-N}M \to 0\qquad \text{and}\quad 0\leq n_k\leq 1=\sup  n_k\]
if on the other side $\gamma'\leq N$ and $q>\frac{N}{\gamma'}$ we have that
\[\|n_k\|_{L^q}=L_k^{-1}\delta_k^{-\frac{N}{q}} \|m_k\|_{L^q} \leq L_k^{\frac{N}{q\gamma'}-1} C_q\to 0  \qquad \text{and}\quad 0\leq n_k\leq 1=\sup  n_k.\]

\noindent {\bf If $\gamma'\leq N$ and $q\leq \frac{N}{\gamma'}$} first of all we observe that, since $\beta<\frac{q}{N}$, 
\[ \|n_k\|_{L^q}=L_k^{-1}\delta_k^{-\frac{N}{q}} \|m_k\|_{L^q} \leq L_k^{\beta\frac{N}{q}-1} C_q\to 0  \qquad \text{and}\quad 0\leq n_k\leq 1=\sup  n_k.\]
We apply \eqref{stimagk} with $r$ as in \eqref{erre} and $s=1$  and we get, using interpolation  between $L^q$ and $L^\infty$ to estimate the norm  $\|m_k\|_{L^r}$, that there holds 
\[\|\tilde g_k\|_\infty\leq \de_k^{\g'}(C_1 \|m_k\|_{L^{\frac{N}{N-2\gamma'} }}+C_2C_q)\leq  L_k^{-\beta \gamma'}(C  L_k^{1-\frac{q}{r}}+C_2C_q)\leq C L_k^{1-\beta \gamma'-\frac{q}{r}}\leq C\]
 since  $\beta\gamma'>1-\frac{q}{r}$. This in turns implies that $|\tilde \lambda_k|\leq C$.

The rest of the proof  follows exactly the same lines of the proof of \cite[Theorem 4.1]{CC}, since we have uniform bounds on $\tilde \lambda_k$ and on $\|\tilde g_k\|_\infty$,  either the $L^1$ or the $L^q$ norm of $n_k$ vanishing as $k\to +\infty$, whereas $\|n_k\|_\infty=1$.  In particular one shows that if $x_k$ is an approximated maximum point of $n_k$ (that is $n_k(x_k)=1-\delta$), then necessarily $\delta_k^{\gamma'+b}|x_k|^b\to +\infty$. If it is not the case, using a priori gradient estimates on $v_k$ as in Proposition \ref{lemma2.5}, we get that $n_k$ is uniformly (in $k$) Holder continuous in the ball $B_1(x_k)$, contradicting the fact that   $n_k\ge0$ and either $\|n_k\|_{L^q}\to 0$ or $\|n_k\|_{L^1}\to 0$. On the other hand, if  $\delta_k^{\g'+b}|x_k|^b\to +\infty$, we may construct a Lyapunov function for the system, which allows for some integral estimates on $n_k$ showing again a uniform (in $k$) Holder bound for $n_k$ in $B_{1/2}(x_k)$ and again getting a contradiction. Therefore one concludes that $L_k\to +\infty$ is not possible. 
\endproof


\section{Pohozaev identity and nonexistence of solutions}

In this section, we study the MFG system \eqref{1} in the case   $V\equiv0$, i.e.
\begin{equation}\label{eps1}
\begin{cases}
-\Delta u+\frac{1}{\g}|\nabla u(x)|^\g+\la=-K_\al\ast m(x)\\
-\Delta m-\mathrm{div} \left(m(x)\nabla u(x)\,|\nabla u(x)|^{\g-2} \right)=0\\
\int_{\R^N}m=M,\quad m\ge 0
\end{cases}\text{in}\,\,\,\R^N.
\end{equation}

\noindent The following Lemma (see Lemma 3.2 in \cite{Cir}) will be useful in order to control the behavior of $m,\,\nabla u,\nabla m$ at infinity.

\begin{lemma}\label{lem3.2}
Let $h\in L^1(\mathbb{R}^N)$. Then, there exists a sequence $R_n\to\infty$ such that 
$$R_n\int_{\partial B_{R_n}}|h(x)|dx\to 0, \quad \text{as}\,\,n\to\infty.$$
\end{lemma}

In order to prove nonexistence of solutions to the MFG system \eqref{eps1} in the \textit{supercritical regime} $0<\al<N-2\g'$, we need a Pohozaev-type identity. 

\begin{prop}[\textbf{Pohozaev identity}]\label{prop3.1} 
Let $(u,m,\la)\in C^2(\R^N)\times W^{1,\frac{2N}{N+\al}}(\R^N)\times \R$ be a solution to \eqref{eps1} such that 
$$m|\nabla u|^\g\in L^1(\R^N) \qquad\text{and}\qquad|\nabla m||\nabla u|\in L^1(\R^N).$$
Then, the following equality holds
\begin{equation}\label{poho}
(2-N)\int\limits_{\R^N}\nabla u\cdot\nabla m \,dx+\left(1-\frac{N} {\g}\right)\int \limits_{\R^N} m|\nabla u|^\g dx=\la N M+\frac{\al+N}{2}\int\limits_{\R^N}m(x)(K_\al \ast m)(x)\,dx.
\end{equation}
\end{prop}

\proof
From Lemma \ref{mC2}, we get that $m$ is twice differentiable, so the following computations are justified. Consider the first equation in \eqref{eps1}, multiplying each term by $\nabla m\cdot x$ and integrating over $B_R(0)$ for $R>0$, we get
\begin{equation}\label{14}
-\int_{B_R}\Delta u\,\nabla m\cdot x\,dx+\frac{1}{\g}\int_{B_R}|\nabla u(x)|^\g \nabla m\cdot x\,dx+\la\int_{B_R} \nabla m\cdot x\,dx=-\int_{B_R}K_\al\ast m(x)\nabla m\cdot x\,dx.
\end{equation}
We take into account each term of \eqref{14} separately. Integrating by parts the first term, we have
\begin{equation}\label{1 di poho}
-\int_{B_R}\Delta u\,\nabla m\cdot x\,dx=\int_{B_R}\nabla u\cdot\nabla(\nabla m\cdot x)dx-\int_{\partial B_R}(\nabla u\cdot \nu)(\nabla m\cdot x)\, d\sigma,
\end{equation}
we observe that
$$\int_{B_R}\nabla u\cdot\nabla(\nabla m\cdot x)dx=\int_{B_R}\sum\limits_{i=1}^N u_{x_i}(\nabla m\cdot x)_{x_i}=\int_{B_R}\nabla u\cdot\nabla m+\int_{B_R}\sum \limits_{i,j}u_{x_i}m_{x_i\,x_j}x_j$$
and integrating by parts the last term of the previous one, we get 
\begin{align*}
\int_{B_R}\sum\limits_{i,j} (u_{x_i}x_j)m_{x_ix_j}&=
\int_{\partial B_R}\sum\limits_{i,j}u_{x_i}m_{x_i}x_j\cdot\frac{x_j}{R}
-\int_{B_R}\sum_{i,j}m_{x_i} u_{x_i\,x_j}x_j-N\int_{B_R}\sum\limits_{i}u_{x_i}m_{x_i}\\
&=\int_{\partial B_R}(\nabla u\cdot\nabla m)x\cdot\nu d\sigma -\int_{B_R}\nabla m\cdot\nabla(\nabla u\cdot x)+(1-N)\int_{B_R}\nabla u\cdot \nabla m.
\end{align*} Note that $x\cdot \nu=R$ on $\partial B_R$. 
Coming back to \eqref{1 di poho} we obtain 
\begin{align}\notag
-\int_{B_R}\Delta u\,\nabla m\cdot x\,dx=&-\int_{B_R}\nabla m\cdot\nabla(\nabla u\cdot x)dx+(2-N)\int_{B_R}\nabla u\cdot \nabla m\,dx \\ \label{uno}
&+\int_{\partial B_R}(\nabla u\cdot\nabla m)(x\cdot\nu) d\sigma -\int_{\partial B_R}(\nabla u\cdot \nu)(\nabla m\cdot x)\,  d\sigma. 
\end{align}
Concerning the second and the third term in \eqref{14}, we get that
\begin{align}\notag
\frac{1}{\g}\int_{B_R}|\nabla u(x)|^\g&\nabla m\cdot x\,dx =\frac{1}{\g}\int_{\partial B_R}|\nabla u|^\g m\, x\cdot\nu \, d\sigma-\frac{1}{\g}\int_{B_R}\,m\,\mathrm{div}(|\nabla u|^\g x)dx=\\\label{3}
&=\frac{1}{\g}\int_{\partial B_R}|\nabla u|^\g m\, x\cdot\nu \, d\sigma-\frac{1}{\g}\int_{B_R}\,m\, \nabla(|\nabla u|^\g)\cdot x\, dx-\frac{N}{\g}\int_{B_R} m|\nabla u(x)|^\g dx
\end{align}
and 
\begin{equation}\label{quattro}
\lambda\int_{B_R}\nabla m \cdot x\,dx=\lambda\int_{\partial B_R}m \,x\cdot\nu d\sigma-\lambda N\int_{B_R} m\, dx.
\end{equation}
Similarly, multiplying the second equation in \eqref{eps1} by $\nabla u\cdot x$ and integrating over $B_R(0)$ we get
\begin{align}\notag
\int_{B_R}\Delta m &\nabla u\cdot x\,dx=-\int_{B_R}\mathrm{div}(m|\nabla u|^{\g-2}\nabla u)\nabla u\cdot x\,dx=\\\notag
&=\int_{B_R}\nabla(\nabla u\cdot x)\cdot(m\,|\nabla u|^{\g-2}\nabla u)\,dx -\int_{\partial B_R}(\nabla u\cdot x)m|\nabla u|^{\g-2}\nabla u\cdot\nu\, d\s=\\\label{6}
&=\int_{B_R}\frac{1}{\g}m\nabla(|\nabla u|^\g)\cdot x\,dx+\int_{B_R} m|\nabla u|^\g dx-\int_{\partial B_R}(\nabla u\cdot x)m|\nabla u|^{\g-2}\nabla u\cdot\nu\, d\s
\end{align}
where we have integrated by parts and then used the following identity
$$\frac{1}{\g}\nabla(|\nabla u|^\g)\cdot x=|\nabla u|^{\g-2}\nabla u\cdot\nabla(\nabla u\cdot x)-|\nabla u|^\g.$$
Integrating by parts the LHS of \eqref{6} we get
\begin{align*}
\int_{\partial B_R}(\nabla m\cdot\nu)&(\nabla u\cdot x)d\s-\int_{B_R} \nabla m\cdot\nabla(\nabla u\cdot x) dx=\\
&=\int_{B_R}\frac{m}{\g}\nabla(|\nabla u|^\g)\cdot x\,dx+\int_{B_R} m|\nabla u|^\g dx- \int_{\partial B_R}(\nabla u\cdot x)m|\nabla u|^{\g-2}\nabla u\cdot\nu\,d\s
\end{align*}
and then isolating the first term in the second line
\begin{align}\notag
-\frac{1}{\g}\int_{B_R}m\nabla(|\nabla u|^\g)\cdot x\,dx=&
\int_{B_R} \nabla m\cdot\nabla(\nabla u\cdot x) dx-\int_{\partial B_R} (\nabla m\cdot\nu)(\nabla u\cdot x)\,d\s\\\label{4}
&+\int_{B_R} m|\nabla u|^\g dx-\int_{\partial B_R}(\nabla u\cdot x)m|\nabla u|^{\g-2}\nabla u\cdot\nu\, d\s
\end{align}
plugging \eqref{4} in \eqref{3} we obtain
\begin{align}\notag
\frac{1}{\g}\int_{B_R}&|\nabla u(x)|^\g\nabla m\cdot x\,dx=\frac{1}{\g}\int_ {\partial B_R}m|\nabla u|^\g x\cdot\nu\,d\s+\int_{B_R}\nabla m\cdot\nabla(\nabla u\cdot x)\,dx+\\\label{due}
&-\int_{\partial B_R} (\nabla m\cdot\nu)(\nabla u\cdot x)\,d\s+\left(1 -\frac{N}{\g}\right)\int_{B_R} m|\nabla u|^\g dx-\int_{\partial B_R}(\nabla u\cdot x)m|\nabla u|^{\g-2}\nabla u\cdot\nu\, d\s.
\end{align}
For what concern the Riesz's potential term, since $m\in L^{\frac{2N}{N+\al}}(\R^N)$ from Theorem \ref{HLS} it follows that $K_\al\ast m\in L^\frac{2N}{N-\al} (\R^N)$, hence by H\"older inequality
$$\bigg|\int_{B_R(0)}K_\al\ast m(x)\nabla m\cdot x\,dx\Big|\le R\int_{B_R(0)}| K_\al\ast m|\,|\nabla m|\,dx\le R\,\|K_\al\ast m\|_{L^\frac{2N}{N-\al}(\R^N)} \,\|\nabla m\|_{L^\frac{2N}{N+\al}(\R^N)},$$
this proves that the term $-\int_{B_R}K_\al\ast m(x)\nabla m\cdot x\,dx$ is finite. We get
\begin{align}\notag
-\int\limits_{B_R}& K_\al\ast m(x)\nabla m\cdot x\,dx=-\int\limits_ {B_R}\int\limits_ {\R^N}\frac{m(y)\nabla m(x)\cdot x}{|x-y|^{N-\al}}dy\,dx=
-\int\limits_{\R^N}\int\limits_{B_R}\frac{m(y)\nabla m(x)\cdot x}{|x-y|^{N-\al}} dx\,dy=\\\label{cinque}
&=-\int\limits_{\R^N}\int\limits_{\partial B_R}\frac{m(x)\,m(y)}{|x-y|^{N-\al}} (x\cdot\nu)\,d\sigma(x)\,dy+\int\limits_{\R^N}\int\limits_{B_R} m(x)\mathrm{div}_x \left(\frac{m(y)}{|x-y|^{N-\al}}x\right)dx\,dy
\end{align}
and furthermore,
\begin{align}\notag
\int\limits_{\R^N}\int\limits_{B_R}&m(x)\,\mathrm{div}_x \left(\frac{m(y)}{|x-y| ^{N-\al}}x\right)dx\,dy=\\\notag
&=(\al-N)\int\limits_{\R^N}\int\limits_{B_R}\frac{m(x)\,m(y)} {|x-y|^{N-\al}}\frac{(x-y) \cdot x}{|x-y|^2}dx\,dy+N\int\limits _{\R^N}\int\limits_ {B_R}\frac{m(x)\,m(y)}{|x-y|^{N-\al}}dx\,dy=\\\label{tre}
&=\frac{\al+N}{2}\int\limits_{\R^N}\int\limits_{B_R}\frac{m(x)\,m(y)}{|x-y|^{N-\al}}dx\,dy+\frac{\al-N}{2}\int\limits_{\R^N}\int\limits_{B_R}\frac{m(x) \,m(y)}{|x-y|^ {N-\al}}\frac{(x+y)\cdot (x-y)}{|x-y|^2}dx\,dy
\end{align}
where we used that $\frac{x\cdot(x-y)}{|x-y|^2}=\frac{1}{2}+\frac{(x+y) \cdot(x-y)}{2|x-y|^2}$. 
Summing up \eqref{uno}, \eqref{quattro}, \eqref{due}, \eqref{cinque} and \eqref{tre} we get the following identity 
\begin{align}\notag
(2-&N)\int\limits_{B_R}\nabla u\cdot\nabla m\,dx+\left(1-\frac{N}{\g}\right)\int \limits_{B_R}m|\nabla u|^\g dx-\la N\int\limits_{B_R}m(x)\,dx\\ \label{pohozaev}
&-\frac{\al+N}{2}\int\limits_{\R^N}\int\limits_{B_R}\frac{m(x)\,m(y)}{|x-y|^{N-\al}}dx\,dy-\frac{\al-N}{2}\int\limits_{\R^N}\int\limits_{B_R}\frac{m(x) \,m(y)}{|x-y|^{N-\al}}\frac{(x+y)\cdot (x-y)}{|x-y|^2}dx\,dy=I_{\partial B_R}
\end{align}
where
\begin{align*}
I_{\partial B_R}=&\int\limits_{\partial B_R}\bigg[-(\nabla u\cdot\nabla m)-\frac{1}{\g}m|\nabla u|^\g-\la\,m\bigg](x\cdot\nu)d\sigma-\int\limits_ {\R^N}\int\limits_{\partial B_R}\frac{m(x)\, m(y)}{|x-y|^{N-\al}} (x\cdot\nu)\,d\sigma(x)\,dy\\
&+\int\limits_{\partial B_R}(\nabla u\cdot\nu)(\nabla m\cdot x)+(\nabla m \cdot\nu)(\nabla u\cdot x)+m|\nabla u|^{\g-2}(\nabla u\cdot x)(\nabla u\cdot \nu)\,d\sigma.
\end{align*}
Now, we let $R$ go to infinity in \eqref{pohozaev}. We observe that (changing variables $x$ and $y$)
$$\int\limits_{\mathbb{R}^N}\int\limits_{\mathbb{R}^N}\frac{m(x) \,m(y)}{|x-y|^{N-\alpha}}\frac{(x+y)\cdot (x-y)}{|x-y|^2}dx\,dy=0.$$
Moreover
$$|I_{\partial B_R}|\le R\int\limits_{\partial B_R}\left(3|\nabla u|\,|\nabla m|+2m\,|\nabla u|^\g+|\la|m\right)d\sigma+\int\limits_{\R^N}R\int\limits_{\partial B_R}\frac{m(x)\,m(y)}{|x-y|^{N-\al}}d\sigma(x)\,dy$$
since by assumption $|\nabla u|\,|\nabla m|$, $m|\nabla u|^\g$ and $m$ $\in L^1(\mathbb{R}^N)$, by Lemma \ref{lem3.2}, we get that for some sequence $R_n\to+\infty$
$$R_n\int\limits_{\partial B_{R_n}}\Big(3|\nabla u|\,|\nabla m|+2m\,|\nabla u|^\gamma+|\lambda|m\Big)d\sigma\to 0,\quad \text{as}\,\,\,n\to+\infty.$$
By means of the same argument, since $m\in L^{\frac{2N}{N+\alpha}}(\R^N)$ implies that $G(x):=\int_{\mathbb{R}^N}\frac{m(x)m(y)}{|x-y|^{N-\alpha}}dy\in L^1(\mathbb{R}^N)$ (by Theorem \ref{HLS2}), we get that there exists a sequence $R_n\to+\infty$ such that
$$R_n\int_{\partial B_{R_n}}G(x)dx \to 0,\quad \text{as}\,\,\,n\to+\infty,$$
which conclude the proof of the Pohozaev-type equality \eqref{poho}.
\endproof

We are now in position to prove nonexistence of classical solutions with prescribed integrability and boundary conditions at $\infty$.

\proof[Proof of Theorem \ref{teo_non_esist}]
We argue by contradiction. Let $(u,m,\la)\in C^2(\R^N)\times W^{1,\frac{2N}{N+\al}} (\R^N)\times\R$ be a solution to \eqref{eps1} such that $u\to +\infty $ as $|x|\to +\infty$  and it holds
$$m|\nabla u|^\g,\,\,\,\,|\nabla m||\nabla u|\,\,\,\in L^1(\R^N).$$
From Proposition \ref{prop3.1} we have the following Pohozaev-type identity
\begin{equation}\label{46}
(2-N)\int\limits_{\R^N}\nabla u\cdot\nabla m \,dx+\left(1-\frac{N} {\g}\right)\int \limits_{\R^N} m|\nabla u|^\g dx=\la N M+\frac{\al+N}{2}\int\limits_{\R^N}m(K_\al\ast m)\,dx.
\end{equation}
Moreover, we obtain the following identities 
\begin{equation}\label{48}
\int\limits_{\R^N}\nabla u\cdot\nabla m\, dx=-\frac{1}{\g}\int\limits_{\R^N}m|\nabla u|^\g dx-\la\,M-\int\limits_{\R^N}\int\limits_{\R^N}\frac{m(x)m(y)}{|x-y|^{N-\al}}dx\,dy
\end{equation}
and
\begin{equation}\label{47}
\int\limits_{\mathbb{R}^N}\nabla u\cdot\nabla m \, dx=-\int\limits_ {\mathbb{R}^N}m|\nabla u|^\gamma\,dx.   
\end{equation}
\textit{Proof of \eqref{48}.} Multiplying the first equation in \eqref{eps1} by $m$ and integrating over $B_R$ we obtain
\begin{equation}\label{48R}
\int\limits_{B_R}\nabla u\cdot\nabla m\,dx-\int\limits_{\partial B_R}m\,\nabla u\cdot\nu\,d\s+\frac{1}{\g}\int\limits_{B_R}m|\nabla u|^\g dx+\la\int\limits_{B_R}m\,dx=-\int\limits_{B_R}\int\limits_{\R^N}\frac{m(x)m(y)}{|x-y|^{N-\al}}dy\,dx.
\end{equation}
By Holder's inequality and using the fact that $m|\nabla u|^\g\in L^1(\R^N)$, we have
$$\int_{\R^N}|\nabla u|m\,dx\le\left(\int_{\R^N}|\nabla u|^\g m\,dx\right)^\frac{1}{\g} \,M^\frac{1}{\g'}<+\infty,$$
hence $|\nabla u|m\in L^1(\R^N)$ and by Lemma \ref{lem3.2}, we get that for some sequence $R_n\to+\infty$
$$\int_{\partial B_{R_n}}m\,\nabla u\cdot \nu\,d\sigma\to 0,\quad \text{as} \,\,\,n\to+\infty.$$
Equality \eqref{48} follows letting $R\to\infty$ in \eqref{48R}. \\
\textit{Proof of \eqref{47}.} For any $s>0$ let us define the set 
$$X_s:=\{x\in\R^N\,|\,u(x)\le s\},$$
and the function 
$$v_s(x):=u(x)-s,\quad\forall x\in\R^N.$$
After a translation we may assume $u(0)=0$. In this way, $\cup_{s>0} X_s=\R^N$, every $X_s$ is non-empty and bounded since $u(x)\to +\infty$ as $|x|\to +\infty$. Multiplying the second equation in \eqref{eps1} by $v_s$ and integrating by parts, we get
$$\int\limits_{X_s}\nabla v_s\cdot\nabla m\, dx=-\int\limits_ {X_s}m\,|\nabla u|^{\gamma-2}\nabla u\cdot \nabla v_s\,dx,$$
since $\nabla v_s=\nabla u$, we obtain \eqref{47} letting $s\to+\infty$.\\
\indent Plugging \eqref{47} in \eqref{48} we get
\begin{equation}
\left(1-\frac{1}{\g}\right)\int\limits_{\R^N}m|\nabla u|^\g dx=\la\,M+\int\limits_{\R^N} \int\limits_{\R^N}\frac{m(x)m(y)}{|x-y|^{N-\al}}dx\,dy   
\end{equation}
and hence
\begin{equation}\label{20}
\int\limits_{\R^N}m|\nabla u|^\g dx=\la\g'M+\g'\int\limits _{\R^N}\int\limits_{\R^N}\frac{m(x)m(y)}{|x-y|^{N-\al}}dx\,dy.
\end{equation}
Using \eqref{47} in \eqref{46}, we have
$$\left(\frac{N}{\g'}-1\right)\int\limits_{\R^N} m|\nabla u|^\g dx=\la N\,M+\frac{\al +N}{2}\int\limits_{\R^N}m(K_\al\ast m)dx$$
and finally from \eqref{20} we obtain
\begin{equation}\label{ug finale}
\left(\frac{N-2\g'-\al}{2}\right)\int_{\R^N}\int_{\R^N}\frac{m(x)m(y)}{|x-y|^{N-\al}}dx\,dy=\g'\la\,M.
\end{equation}
Recall that by Lemma \ref{lem3.3}, we have that $\la\le0$ and by assumption $N-2\g'-\al>0$, so we get a contradiction.
\endproof

\begin{obs}
One could observe that the the previous proof (with slight changes) holds
also in the case when $u\to-\infty$ as $|x|\to+\infty$, hence one may ask why we do not consider this possibility. This is due to the fact that the property of ergodicity for the process is strictly related to the existence of a Lyapunov function (refer to \cite{H}). More in detail, a necessary condition to have an ergodic process is  
$$\nabla u\cdot x>0, \quad\text{for $x$ large}$$
(see also \cite{Ci} and references therein). As a consequence, the case $u\to-\infty$ as $|x|\to+\infty$ is not relevant.
\end{obs}
%


\section{Existence of classical solutions to the MFG system} \label{sez pto fix}

First of all we consider a regularised version of problem \eqref{1}, namely

\begin{equation}\label{1reg}
\begin{cases}
-\Delta u+\frac{1}{\g}|\nabla u|^\g+\la=V(x)-K_\al\ast m\ast\p_k (x)\\
-\Delta m-\mathrm{div}\left(m(x)\nabla u(x)\,|\nabla u(x)|^{\g-2} \right)=0\\ \int_{\R^N}m=M,\quad m\ge 0
\end{cases}\quad\text{in}\quad\R^N
\end{equation}

\noindent where $\p_k$ is a sequence of standard symmetric mollifiers approximating the unit as $k\to+\infty$ (i.e. a sequence of symmetric functions on $\R^N$ such that $\p_k\in C^\infty_0(\R^N)$, $\mathrm{supp}\,\p_k\subset\overline{B_{1/k}(0)}$, $\int\p_k=1$ and $\p_k\ge0$). For every $k$ fixed, using  Schauder Fixed Point Theorem, we will prove existence of $(u_k,m_k,\la_k)$ solution to \eqref{1reg},    and then, exploiting a priori uniform estimates on these solutions,   we will show that we may pass to the limit as  $k\to+\infty$ and get  a solution of the MFG system \eqref{1}.

\subsection{Solution of the regularized problem}

We consider \eqref{1reg} with $k$ fixed: 
\begin{equation}\label{1 phi}
\begin{cases}
-\Delta u+\frac{1}{\g}|\nabla u|^\g+\la=V(x)-K_\al\ast m\ast\p(x)\\
-\Delta m-\mathrm{div}\left(m(x)\nabla u(x)\,|\nabla u(x)|^{\g-2} \right)=0\\ \int_{\R^N}m=M,\quad m\ge 0
\end{cases}\quad\text{in}\quad\R^N
\end{equation}
We are going to construct solution to \eqref{1 phi} by using the following
version of the well-known Schauder Fixed Point Theorem. Construction of solutions to MFG systems by exploiting fixed point arguments is quite classical in the literature, see \cite{BaFe, Cir, GoPibook, LL1}.
\begin{theorem}[Corollary 11.2 \cite{GT}]\label{cor11.2GT}
Let $A$ be a closed and convex set in a Banach space $X$ and let $\mathcal{F}$ be a continuous map from $A$ into itself such that the image $\mathcal{F}(A)$ is precompact. Then, $\mathcal{F}$ has a fixed point.
\end{theorem}

Let $\xi,C>0$ (which will be chosen later), $M>0$ and $\bar{p}>\frac{N}{\al}$, we define the set
\begin{equation}\label{setA} 
A_{\xi,M,C}:=\left\{\mu\in L^{\bar{p}}(\R^N)\cap L^1(\R^N)\,\bigg|\,\|\mu\|_ {L^\frac{2N}{N+\al}(\R^N)}\le\xi,\,\,\,\int_{\R^N}\mu\,dx=M,\,\,\,\mu\ge0,\,\,\,\int_{\R^N}\mu V(x)\,dx \le C\right\}.
\end{equation} 

\begin{lemma}\label{A chiuso}
For any choice of $\xi,M,C>0$, the set $A_{\xi,M,C}\subset L^{\bar{p}}(\R^N)$ is closed and convex.
\end{lemma}

\proof
The set $A_{\xi,M,C}$ is convex since it is intersection of convex sets.\\
\indent Let now  $(\mu_n)_n$ be a sequence in $A_{\xi,M,C}$ which converges in $L^{\bar{p}}$ to   $\bar{\mu}$. Obviously $\bar{\mu}\ge0$ and since $\mu_n\rightharpoonup\bar{\mu}$ in $L^\frac{2N}{N+\al}(\R^N)$ by weak lower semicontinuity of the norm we have that 
$$\|\bar{\mu}\|_{L^{\frac{2N}{N+\al}}(\R^N)}\le \mathrm{lim\,inf}\,\|\mu_n\|_{L^\frac{2N} {N+\al}(\R^N)}\le\xi.$$
From Fatou's Lemma we get also that $\int_{\R^N}\bar{\mu}V(x)\le\mathrm{lim\,inf}\int_{\R^N} \mu_n V(x)\le C$. 
Note that due to the fact that $0\leq \int_{\R^N} \mu_n V(x)\le C$, and that $V$ is coercive, see \eqref{V}, $\mu_n$ are uniformly integrable, since for every $R>>1$,
$0\leq   \int_{|x|\geq R} \mu_ndx \leq \frac{C_V}{R^b}\int_{\R^N} V(x)\mu_n dx\leq \frac{CC_V}{R^b}$. Due to the fact that $\mu_n\to \bar \mu$ in $L^{\bar p}$, we have also that they have   uniformly absolutely continuous integrals, so we may apply  the Vitali convergence theorem and obtain that  $\mu_n\to\bar{\mu}$ in $L^1(\R^N)$ and hence $\int_{\R^N}\bar{\mu}\,dx=M$. This proves that  $\bar{\mu}\in A_{\xi,M,C}$, and hence that $A_{\xi,M,C}$ is closed. 
\endproof

We define the map $F:A_{\xi,M,C}\to C^2(\R^N)\times\R$ which to every element $\mu\in A_{\xi,M,C}$ associates a solution $(u,\bar{\la})\in C^2(\R^N)\times\R$ to the HJB equation
\begin{equation}\label{eq1_Schauder}
-\Delta u+\frac{1}{\g}|\nabla u|^\g+\la=V(x)-K_\al\ast\mu\ast\p(x),\quad \text{in}\,\,\,\R^N
\end{equation}
where $\bar{\la}$ is defined as in \eqref{lambda bar} (refer to \cite{BM}); and the map $G$ which to the couple $(u,\bar{\la})$ associates the function $m$ which solves (weakly)
\begin{equation}\label{eq2_Schauder}
\begin{cases}-\Delta m-\mathrm{div}\left(m(x)\nabla u(x)\,|\nabla u(x)|^{\g-2} \right)=0\\ \int_{\R^N}m=M,\quad m\ge 0\end{cases}.
\end{equation}
We look for a fixed point of the map $\mathcal{F}:\mu\mapsto m$ defined as the composition of $F$ and $G$, namely $\mathcal{F}(\mu):=G(F(\mu))$.

We are going to show that, once we have fixed $M$ (in an appropriate range), it is possible  to choose appropriately $\xi$ and $C$ in such a way that the map $\mathcal{F}$  defined on $A_{\xi, M, C}$ satisfies the assumptions of the Schauder Fixed Point Theorem \ref{cor11.2GT}. As we will see the regularization with $\p$ in the system \eqref{1 phi} is necessary in order to get precompactness of the image of $\mathcal{F}$. 
We start with  some preliminary results.

\begin{prop}\label{prop fix}
Let us consider $\mu\in A_{\xi,M,C}$,  $(u,\bar{\la})=F(\mu)$ and $m=G(u, \bar \lambda)=\mathcal{F}(\mu)$. Then,
\begin{itemize}
    \item[i)] there exists a constant $C>0$ such that
          \begin{equation}\label{u1b}
          |\nabla u(x)|\le C(1+|x|)^\frac{b}{\g}
          \end{equation}
          where $C$ depends on $C_V,b,\g,N,\bar{\la},\|K_\al\ast\mu\ast\p\|_\infty$.
    \item[ii)] the function $u$ is unique up to addition of constants and there
          exists $C>0$ such that
          \begin{equation}\label{u2b}
          u(x)\ge C|x|^{\frac{b}{\g}+1}-C^{-1}.
          \end{equation}
    \item[iii)] it holds
          \begin{equation}\label{stima bar lambda basso alto}
          -K_1\le\bar{\la}\le K_2
          \end{equation}
          where $K_1$ and $K_2$ are positive constants depending respectively on $\|K_\al\ast\mu\ast\p\|_\infty$ and on $C_V,b,\g,N$.
   \item[iv)] the function $m$ is unique, $m\in(W^{1,1}\cap L^\infty)(\R^N)$, 
          $\sqrt{m}\in W^{1,2}(\R^N)$, $m\in W^{1,p}(\R^N)$ $\forall p>1$ and it holds
          \begin{equation}\label{stima nabla m fix}
          \|\nabla m\|_{L^p(\R^N)}\le C\|m^{\frac{1}{p}}|\nabla u|^{\g-1} \|_{L^p(\R^N)} \|m^{1-\frac{1}{p}}\|_{L^\infty(\R^N)}.   
          \end{equation}
          Moreover, the following integrability properties are verified
          \begin{equation}\label{u3b}
          m|\nabla u|^\g\in L^1(\R^N),\qquad m V\in L^1(\R^N),\qquad |\nabla u|\,|\nabla m|\in L^1(\R^N).
          \end{equation}
\end{itemize}
\end{prop}

\proof
\indent $i)$ Since $\mu\ast\p\in L^1(\R^N)\cap L^{\bar p}(\R^N)$ with $\bar p>\frac{N}{\alpha}$, by  Theorem \ref{holderRiesz}  we obtain  $K_\al\ast(\mu\ast\p)\in C^{0,\theta}(\R^N)$ for some $\theta\in (0,1)$ and   $\|K_\al\ast \mu\ast\p\|_\infty\le C_{N,\al,\bar{p}}\|\mu\ast\p\|_{L^{\bar{p}}}+\|\mu\ast\p\|_{L^1}\le C_{N,\al,\bar{p}} \|\mu\|_{L^{\bar{p}}}+M$.  We can apply therefore  Proposition \ref{lemma2.5}, which gives us the following estimate
\begin{equation}\label{stime u}
|\nabla u(x)|\le C(1+|x|)^\frac{b}{\g}
\end{equation}
where $C$ is a constant depending on $C_V, b, \g, N, \bar{\la}, \|K_\al\ast(\mu\ast\p)\| _\infty$. This proves \eqref{u1b}.\\
\indent $ii)$ Since, by construction, $u$ is a solution to \eqref{eq1_Schauder} with $\la=\bar{\la}$ then by Proposition \ref{teo2.7} $ii)$ it follows uniqueness up to additive constants and \eqref{u2b}.\\
\indent $iii)$ The fact that $\bar \lambda\leq  K_2$ is a direct consequence of  Lemma \ref{lem3.3}. Furthermore, if $\bar{x}$ is a minimum point of $u$, evaluating \eqref{eq1_Schauder} at $\bar{x}$ we have that
$$\bar{\la}\ge V(\bar{x})-K_\al\ast\mu\ast\p(\bar{x})\ge-\|K_\al\ast\mu\ast\p\|_{\infty}\ge-K_1$$
since $V(x)\ge0$ in $\R^N$.\\
\indent $iv)$ For $r>1$, let us consider the function $h(x):=u(x)^r$, one can observe that 
\begin{align*}
-\Delta h+|\nabla& u|^{\g-2}\nabla u\cdot\nabla h= ru^{r-1}\left(-(r-1)\frac{|\nabla u|^2}{u}-\Delta u+|\nabla u|^\g\right)\\
&=ru^{r-1}\left(-\Delta u+\frac{1}{\g}|\nabla u|^\g-(r-1)\frac{|\nabla u|^2}{u}+\frac{1}{\g'}|\nabla u |^\g\right)\\
&=ru^{r-1}\left(-\bar{\la}+V-K_\al\ast\mu\ast\p-(r-1)\frac{|\nabla u|^2}{u}+\frac{1}{\g'}|\nabla u |^\g\right),
\end{align*}
where in the last equality we used the fact that $u$ solves \eqref{eq1_Schauder}. 
Denoting by 
$$H(x):=-\bar{\la}+V(x)-K_\al\ast\mu\ast\p(x)-(r-1)\frac{|\nabla u|^2}{u}+\frac{1}{\g'}|\nabla u |^\g,$$
from \eqref{stima bar lambda basso alto}, \eqref{V} and the fact that $K_\al\ast\mu\ast\p \in L^\infty$, we get 
$$H(x)\ge (r-1)|\nabla u|^\g\left(\frac{1}{\g'(r-1)}-\frac{|\nabla u|^{2-\g}}{u}\right)+C_V^{-1}|x|^b-C\ge1,\quad\text{for}\,\,\,|x|>R$$
taking $R$ sufficiently large. Hence, for $|x|>R$
$$-\Delta h+|\nabla u|^{\g-2}\nabla u\cdot\nabla h\ge C|x|^{(\frac{b}{\g}+1)(r-1)}>0$$
this means that $h$ is a Lyapunov function for the stochastic process with drift $|\nabla u|^{\g-2}\nabla u$. Since $m$ solves \eqref{eq2_Schauder}, it is the density of the invariant measure associated to this process. So, from \cite[Proposition 2.3]{MPR} we get that 
\begin{equation}\label{C|x|^b/g}
m |x|^{(\frac{b}{\g}+1)(r-1)}\in L^1(\R^N).
\end{equation}
More in general, since the value of $r>1$ can be chosen arbitrarily, from \eqref{C|x|^b/g} we have that for any $q>0$ 
\begin{equation}\label{m|x|^qL1}
m|x|^q\in L^1(\R^N),
\end{equation}
in particular $m|x|^b\in L^1(\R^N)$, so taking into account estimates \eqref{stime u} and \eqref{V} we obtain that 
$$m|\nabla u|^\g\in L^1(\R^N) \qquad \text{and} \qquad mV\in L^1(\R^N).$$
With the same argument (since $|\nabla u|^{p(\g-1)}$ has polynomial growth) it follows that 
$$m|\nabla u|^{p(\g-1)}\in L^p(\R^N),\quad\forall p>1$$
hence from \cite[Corollary 3.2 and Theorem 3.5]{MPR} we get that 
$$m\in W^{1,1}(\R^N)\cap L^\infty(\R^N).$$
Moreover, using the fact that $m$ is a weak solution to the Kolmogorov equation and  H\"older inequality, we obtain that for any $\phi\in C^\infty_0(\R^N)$ we have
$$\bigg|\int_{\R^N}m\Delta\phi\,dx\bigg|\le\int_{\R^N}m|\nabla u|^{\gamma-1}|\,|\nabla\phi|\,dx\le\|m^{\frac{1}{p}}|\nabla u|^{\g-1}\|_p\|m^{1-\frac{1}{p}}\|_\infty\| \nabla\phi\|_{p'}.$$ 
Since $m^\frac{1}{p}|\nabla u|^{\g-1}\in L^p(\R^N)$ and $m^{1-\frac{1}{p}}\in L^\infty(\R^N)$, by Proposition \ref{prop2.4CC} we get that 
$$m\in W^{1,p}(\R^N), \,\,\,\forall p>1$$
and estimate \eqref{stima nabla m fix} holds. Finally, from \cite[Theorem 3.1]{MPR} we have that $\sqrt{m}\in W^{1,2}(\R^N)$ and $\int_{\R^N}\frac{|\nabla m|^2}{m}<+\infty$, so using H\"older inequality we obtain
$$\int_{\R^N}|\nabla u|\,|\nabla m|\le\big\|\,|\nabla u|\sqrt{m}\,\big\|_2 \,\bigg\| \frac{|\nabla m|}{\sqrt{m}}\bigg\|_2<+\infty.$$
Since the function $u$ is unique up to additive constants, $\nabla u$ is fixed and hence, by existence of a Lyapunov function, it follows immediately uniqueness of $m$ solution to the Kolmogorov equation.
\endproof


We show now that once we fix the mass $M$ (in $(0, +\infty)$ in the mass-subcritical case, or below a certain threshold in the mass-supercritical and mass-critical regime), then we may choose the constant $\xi, C$ in the definition \eqref{setA} of the set $A_{\xi,M,C}$ such that the map $\mathcal{F}$ maps $A_{\xi,M,C}$ into itself. 

\begin{lemma}\label{K in se}
We have the following results:
\begin{itemize}
    \item[i)] if $N-\g'<\al<N$, then for any $M>0$, there exist $\xi, C>0$ such that  $\mathcal{F}$ maps $A_{\xi,M,C}$ into itself;
    \item[ii)] if $N-2\g'\le\al\leq N-\g'$ then there exists a positive real value $M_0=M_0(N,\al,\g, C_V, b)$ such that if $M\in\left(0,M_0\right)$ there exist $\xi, C>0$ such that  $\mathcal{F}$ maps the set $A_{\xi,M,C}$ into itself.\end{itemize}
\end{lemma}

\proof
Let $\mu\in A_{\xi,M,C}$, $m=\mathcal{F}(\mu)$ and $(u,\bar{\la})=F(\mu)$ as above.
Since by Proposition \ref{prop fix} $iv)$ $m\in L^\infty(\R^N)$, by interpolation it follows that $m\in L^{\bar{p}}(\R^N)$.
Multiplying \eqref{eq1_Schauder} by $m$ and integrating over $B_R$, we obtain
$$-\int_{B_R}m\Delta u\,dx+\frac{1}{\g}\int_{B_R}m|\nabla u|^\g dx+\bar{\la}\int_{B_R} m\,dx=\int_{B_R}V(x)m\,dx-\int_{B_R}m(K_\al\ast \mu\ast\p)\,dx$$
and integrating by parts the first term
\begin{align}\notag
\int_{B_R}\nabla m\cdot\nabla u\,dx-\int_{\partial B_R}m\nabla u\cdot\nu\,d\s+ \frac{1}{\g}\int_{B_R}&m|\nabla u|^\g dx+\bar{\la}\int_{B_R} m\,dx\\ \label{53 integr}
&=\int_{B_R}V(x)m\,dx -\int_{B_R}m(K_\al\ast\mu\ast\p)\,dx.
\end{align}
From the fact that $\int_{\R^N}m=M$ and $m|\nabla u|^\g\in L^1(\R^N)$, by H\"older inequality we get that $m|\nabla u|\in L^1(\R^N)$, hence by Lemma \ref{lem3.2} for some sequence $R_n\to+\infty$ we have that $\int_{\partial B_{R_n}}m\nabla u\cdot\nu\,d\sigma \to 0$. Since $m(K_\al\ast\mu\ast\p)\in L^1(\R^N)$ and \eqref{u3b} holds, letting $R$ go to $+\infty$ in \eqref{53 integr} we obtain that
\begin{equation}\label{eq1}
\int_{\R^N}\nabla u\cdot\nabla m\,dx=-\frac{1}{\g}\int_{\R^N}m|\nabla u|^\g dx-\bar{\la}M+\int_{\R^N}V(x)m\,dx-\int_{\R^N}m(K_\al\ast\mu\ast\p)\,dx.
\end{equation}
Moreover, from the fact that $m$ solves (weakly) the Kolmogorov equation in \eqref{eq2_Schauder}, following the proof of identity \eqref{47}, we have that
\begin{equation}\label{eq2}
\int_{\R^N}\nabla u\cdot\nabla m\,dx=-\int_{\R^N}m|\nabla u|^\g dx.
\end{equation}
Putting together \eqref{eq1} and \eqref{eq2} we get that
\begin{equation}\label{=}
\frac{1}{\g'}\int_{\R^N}m|\nabla u|^\g dx+\int_{\R^N}mV\,dx=\bar{\la}M+\int_{\R^N}m (K_\al\ast\mu\ast\p)\,dx.
\end{equation}
Since  $\bar{\la}\le K_2$ (from \eqref{stima bar lambda basso alto}), using \eqref{disug hls} we have
\begin{align}\notag
\int_{\R^N}m|\nabla u|^\g dx &\le C_1 M+\,C_2 \,\|m\|_{L^\frac{2N}{N+\al}(\R^N)} \|\mu\ast\p\|_{L^\frac{2N}{N+\al}(\R^N)}\\ \notag
&\le C_1 M+\,C_2 \,\|m\|_{L^\frac{2N}{N+\al}(\R^N)}\|\mu\|_{L^\frac{2N}{N+\al}(\R^N)}\\\label{stima alto}
&\le C_1 M+\,C_2\,\xi \,\|m\|_{L^\frac{2N}{N+\al}(\R^N)}
\end{align}
where $C_1=C_1(\g,C_V,b,N)$ and $C_2=C_2(\al,N,\g)$.\\
\noindent {\bf Choice of $\xi$}. 
First of all we show that we may choose $\xi$ in such a way that if $\mu\in A_{\xi, M,C}$ then $\|m\|_{L^{\frac{2N}{N+\al}}(\R^N)}=\|\mathcal{F}(\mu)\|_{L^{\frac{2N}{N+\al}}(\R^N)}\leq \xi$.\\
\indent Let us fix $a:=\frac{2\g'}{N-\al}$. Notice that $a>2$ if $\al>N-\g'$, $a=2$ if $\al=N-\g'$, $a\in(1,2)$ if $N-2\g'<\al<N-\g'$ and $a=1$ when $\al=N-2\g'$.

In the case when $N-\g'\le\al<N$, using estimate \eqref{stima m_p 1}, we get
\begin{equation}\label{stima m^a 1}
\|m\|^a_{L^\frac{2N}{N+\al}(\R^N)}\le CM^{a-1}\int_{\R^N}m|\nabla u|^\g dx  
\end{equation}
where $C$ is a constant depending on $N$, $\al$ and $\g$; whereas if $N-2\g'\le\al<N-\g'$ using estimate \eqref{stima m_p 2}, we get
\begin{equation}\label{stima m^a 2}
\|m\|^a_{L^\frac{2N}{N+\al}(\R^N)}\le CM^{a-1}\left(\int_{\R^N}m|\nabla u|^\g dx+M\right)
\end{equation}
where $C$ is a constant depending on $N$, $\al$ and $\g$.

From \eqref{stima alto} and  either \eqref{stima m^a 1} or \eqref{stima m^a 2} we obtain that 
\begin{equation}\label{71}
\|m\|^a_{L^\frac{2N}{N+\al}(\R^N)}\le C_1 M^a+C_2 M^{a-1}\xi\,\|m\|_ {L^\frac{2N}{N+\al} (\R^N)}.   
\end{equation}
\indent We define the function 
\[f(t):= t^a-C_2 M^{a-1}\xi t-C_1 M^a \] and observe that \eqref{71} is equivalent to $f(\|m\|_{L^\frac{2N}{N+\al}(\R^N)})\leq 0$. 
When $a>1$, $f(\|m\|_{L^\frac{2N}{N+\al}(\R^N)})\leq 0$ is equivalent to $\|m\|_{L^\frac{2N}{N+\al}(\R^N)}\leq t_0$, where $t_0$ is the unique zero of $f$. So, in order to  conclude that $ \|m\|_{L^\frac{2N}{N+\al}(\R^N)}\leq \xi$ it is sufficient to choose $\xi$ such that $f(\xi)\geq 0$.\\
\indent\textit{Case $N-\g'<\al<N$.} In this case 
since $a>2$ and  $f(\xi)= \xi^a-C_2 M^{a-1}\xi^2 -C_1 M^a$ then  for every fixed $M>0$, there exists  $\xi_M$ such that $f(\xi)\geq 0$ for every $\xi\geq \xi_M$ and we have done.\\
\indent\textit{Case $\al=N-\g'$.} In this case $a=2$, so arguing as before, and recalling that $f(\xi)=\xi^2-C_2M\xi^2-C_1M^2$,  we get that whenever $M<  \frac{1}{C_2}:=M_0$ there exists $\xi_M$ such that $f(\xi)\geq 0$ for every $\xi\geq \xi_M$.\\
\indent\textit{Case $N-2\g'<\al<N-\g'$.} In this case $a\in (1,2)$. Denote $g(t):=t^a-C_2 M^{a-1}t^2-C_1 M^a$. We aim to find $\xi$ such that $g(\xi)\geq 0$. This is possible if and only if $g(t_{max})\geq 0$, where $t_{max}=\left(\frac{a}{2C_2 M^{a-1}}\right)^\frac{1}{2-a}$ is the maximum point of $g$.  Evaluating $g$ in this point we get
$$g(t_{max})=\left(\frac{2-a}{2}\right)C_2^{-\frac{a}{2-a}}\left(\frac{a}{2}\right)^\frac{a}{2-a}M^\frac{a(1-a)}{2-a}-C_1M^a.$$
Since $a\in(1,2)$ we have that $\frac{2-a}{2}>0$ and $\frac{a(1-a)}{2-a}<0$, hence 
$$g(t_{max})\ge0$$
provided that $M$ is sufficiently small. We may choose $\xi=t_{max}$ (or more generally $\xi$ in the range of values $t$ such that $g(t)\ge0$).\\ 
\indent\textit{Case $\al=N-2\g'$.} Since $a=1$,  the function $f$ reads $f(t)=t-C_2 \xi t-C_1 M$ and since  $f(\|m\|_{L^\frac{2N}{N+\al}(\R^N)})\leq 0$ we have, if $\xi<\frac{1}{C_2}$, 
\begin{equation}\label{a=1}
\|m\|_{L^\frac{2N}{N+\al}(\R^N)}\le \frac{C_1 M}{1-C_2\xi}.
\end{equation}
We look for some condition on $M$ under which we may choose $\xi$ such that  $\frac{C_1 M}{1-C_2\xi}\le\xi$. Observe that this is equivalent to  $C_2\xi^2-\xi+C_1 M\leq 0$. If $M\le\frac{1}{4C_1C_2}$, then it is sufficient to choose $\xi$ in the range $\left[\frac{1-\sqrt{1-4C_1C_2 M}}{2C_2},\frac{1+\sqrt{1-4C_1C_2 M}}{2C_2}\right]\cap \left(0, \frac{1}{C_2}\right).$ \\
\noindent {\bf Choice of $C$}. 
Notice that in each of the previous cases, from \eqref{=}, \eqref{stima alto} and the fact that $\|m\|_{L^\frac{2N}{N+\al}(\R^N)}\le\xi$, we get
\begin{align}\notag
\int_{\R^N}m V dx \le C_1 M+C_2\xi^2.
\end{align}
So it is sufficient to choose $C$ greater or equal to $C_1 M+C_2\xi^2$.\\
\noindent We can conclude that $\mathcal{F}$ maps the set $A_{\xi,M,C}$ into itself.
\endproof
 
We show now that the image of $\mathcal{F}$ is precompact, that is  relatively compact. Here is the main point in which the regularization with the mollifier $\p$ comes into play. 

\begin{lemma}\label{pre comp}
Let  $M$ and $\xi, C$ as given by Lemma \ref{K in se}. Then the image $\mathcal{F}({A_{\xi,M,C}})$ is precompact.
\end{lemma}

\proof
Let us consider a sequence $(m_n)_n\subset\mathcal{F}(A_{\xi,M,C})$, in order to prove that $\mathcal{F}(A_{\xi,M,C})$ is precompact in $A_{\xi,M,C}$, we have to show that that $(m_n)_n$ admits a subsequence converging in $L^{\bar{p}}$-norm to a point belonging to $A_{\xi,M,C}$. There exists a sequence $(\mu_n)_n\subset A_{\xi,M,C}$ such that $\mathcal{F}(\mu_n)=m_n$ for every $n\in \mathbb{N}$, considering also $(u_n,\bar{\la}_n)=F(\mu_n)$, we have that for every $n\in\mathbb{N}$ the triple $(u_n,m_n,\bar{\la}_n)$ is such that
$$\begin{cases}
-\Delta u_n+\frac{1}{\g}|\nabla u_n|^\g+\bar{\la}_n=V(x)-K_\al\ast\mu_n\ast\p(x)\\
-\Delta m_n-\mathrm{div}(m_n\nabla u_n|\nabla u_n|^{\g-2})=0\\
\int_{\R^N}m_n=M \quad m_n\ge0. 
\end{cases}.$$
Note that by Young's convolution inequality $\|\mu_n\ast\p\|_{L^{q}(\R^N)}\leq \|\mu_n\|_{L^{1}(\R^N)} \|\p\|_{L^{q}(\R^N)} = M \|\p\|_{L^{q}(\R^N)} $ for every $q$. Therefore by Proposition \ref{holderRiesz} we get that $K_\al\ast\mu_n\ast\p\in L^\infty\cap C^{0, \theta}$ for some $\theta\in (0,1)$ uniformly in $n$, that is $\|K_\al\ast\mu_n\ast\p\|_{L^\infty(\R^N)}\leq C$, for some $C$ independent of $n$. By Proposition \ref{prop fix} we have that $u_n$ are bounded from below, that $m_n\in L^\infty$ and that $\bar{\la}_n$ are equibounded in $n$, so we may apply Theorem \ref{teo4.1new} (actually a simpler version, with $W_n(x)=V(x)-K_\al\ast\mu_n\ast\p$ and $G_{k,\al}\equiv 0$). So we  obtain that there exists a positive constant $C$ not depending on $n$ such that
\begin{equation}\label{m n infty<C}
\|m_n\|_{L^\infty(\R^N)}\le C, \qquad\forall n\in\mathbb{N}.
\end{equation}
Now we use  Proposition \ref{kolm} $ii)$, since $m_n\in L^q(\R^N)$ (where $q$ is defined as in Proposition \ref{kolm}) and $E_n\le C_1M+C_2\xi^2$ and we get that 
$$\|m_n\|_{W^{1,\ell}(\R^N)}\le C, \quad \forall \ell<q$$
where the constant $C$ does not depend on $n$. Hence, by Sobolev compact embeddings, $m_n\to\bar{m}$ strongly in $L^s(K)$ for $1\le s<q^*$ and for every $K\subset\subset\R^N$. Moreover, using the fact that $\int_{\R^N}m_n V\,dx\le C$ uniformly in $n$ and \eqref{V} we get that for $R>1$
$$C\ge\int_{\R^N}m_nV\,dx\ge\int_{|x|\geq R}m_nV\,dx\ge CR^b\int_{|x|\geq R}m_n(x)\,dx$$
that is 
$$\int_{|x|\geq R}m_n(x) dx\to 0,\,\,\,\text{as}\,\,\,R\to+\infty.$$
Using also the uniform estimate \eqref{m n infty<C}, from the Vitali Convergence Theorem we obtain that up to sub-sequences
\begin{equation}\label{L1}
m_n\to\bar{m}\quad\text{in} \,\,\, L^1(\R^N)
\end{equation}
and as a consequence $\int_{\R^N}\bar{m}(x)dx=M$. Finally, from \eqref{m n infty<C} and \eqref{L1}, we deduce that $m_n\to\bar{m}$ strongly in $L^{\bar{p}}(\R^N)$. 
Since $A_{\xi,M,C}$ is closed and by Lemma \ref{K in se} we have that $\mathcal{F}(A_{\xi,M,C})\subset A_{\xi,M,C}$, we may conclude that $\bar{m}\in A_{\xi,M,C}$.
\endproof

Finally we show that $\mathcal{F}$ is continuous. 

\begin{lemma}\label{F continua}
Let $\xi$, $M$ and $C$ as given by Lemma \ref{K in se}. Then, the map $\mathcal{F}$ is continuous.
\end{lemma}

\proof
Let $(\mu_n)_n$ be a sequence in $A_{\xi,M,C}$ such that $\mu_n\to\tilde{\mu}\in A_{\xi,M,C}$ strongly in $L^{\bar{p}}(\R^N)$. In order to prove that the map $\mathcal{F}$ is continuous, we have to show that $\mathcal{F}(\mu_n)\to\mathcal{F}(\tilde{\mu})$ with respect to the $L^{\bar{p}}$-norm, that is $m_n\to\tilde{m}$ strongly in $L^{\bar{p}}(\R^N)$.\\
We consider the sequence made by the couples $(u_n,\bar{\la}_n)\in C^2(\R^N)\times\R$, where $(u_n,\bar{\la}_n)=F(\mu_n)$ $\forall n\in\mathbb{N}$, as previously defined.
As observed in  Lemma \ref{pre comp}, $K_\al\ast\mu_n\ast\p$ is uniformly bounded in $L^\infty$. So by Proposition \ref{prop fix} we have that $\bar\la_n$ are uniformly bounded, that
$$|\nabla u_n(x)|\le C(1+|x|^\frac{b}{\g})\qquad\text{uniformly in }n$$
and then consequently
\begin{equation}\label{delta u_n}
|\Delta u_n|\le C(1+|x|^b),\qquad\text{uniformly in }n.
\end{equation}
Up to extracting a subsequence we assume that $\bar\la_n\to\la^{(1)}$. 
Since $u_n$ is a classical solution to the HJB equation, by classical elliptic regularity estimates applied to $v_n(x):=u_n(x)-u_n(0)$ (refer e.g to \cite[Theorem 8.32]{GT}) for any $\theta\in(0,1]$ and $K\subset\subset\R^N$ we get
$$\|v_n\|_{C^{1,\theta}_{loc}(K)}\le C\qquad\text{uniformly in }n$$
(notice that the previous estimate holds for $\theta=1$ thanks to \eqref{delta u_n}).
By Arzelà-Ascoli Theorem, up to extracting a subsequence, we get that
$$v_n\to u^{(1)}\qquad \text{locally uniformly in }C^{1,\theta}$$
and hence
$$\nabla u_n\to \nabla u^{(1)}\qquad \text{locally uniformly in }C^{0,\theta}.$$
Since $\|(\mu_n-\tilde \mu)\ast\p\|_{L^{\bar p}(\R^N)}  \leq 2 M   \|\p\|_{L^{\bar p}(\R^N)} $, by  Theorem \ref{holderRiesz}   we get that
$$\|K_\al\ast\p\ast\mu_n\|_{C^{0,\al-N/\bar{p}}}\le C,\qquad\text{uniformly in }n$$
and 
$$\|K_\al\ast\p\ast\mu_n-K_\al\ast\p\ast\tilde{\mu}\|_{L^\infty(\R^N)}\le C_{N,\al,\bar{p}}\|\mu_n-\tilde{\mu}\|_{L^{\bar{p}}(\R^N)}+\|\mu_n-\tilde{\mu}\|_{L^1(\R^N)}.$$
Since $\mu_n\to\tilde{\mu}$ in $L^1(\R^N)\cap L^{\bar{p}}(\R^N)$, then up to subsequences
$$K_\al\ast\p\ast\mu_n\longrightarrow K_\al\ast\p\ast\tilde{\mu}\quad 
\text{locally uniformly in $\R^N$}.$$
By stability   with respect to locally uniform convergence, we get that $(u^{(1)},\la^{(1)})$ is a solution (in the viscosity sense) to the HJB equation 
$$-\Delta u+\frac{1}{\g}|\nabla u|^\g+\la=V(x)-K_\al\ast\p\ast\tilde{\mu}(x), \quad\text{on} \,\, \R^N.$$
Let $(\tilde{u},\tilde{\la})=F(\tilde{\mu})$, we want to show that $\tilde{\la}=\la^{(1)}$. Assume by contradiction that $\tilde{\la}\not=\la^{(1)}$, without loss of generality we can assume that $\la^{(1)}<\tilde{\la}-2\e$ for a certain $\e>0$. Then, for $n$ sufficiently large $\bar{\la}_n<\tilde{\la}-\e$ and, possibly enlarging $n$, we have also $\|K_\al\ast\p\ast\mu_n-K_\al\ast\p \ast\tilde{\mu}\|_\infty\le\e$. One can observe that 
$$-\Delta\tilde{u}+\frac{1}{\g}|\nabla\tilde{u}|^\g+\tilde{\la}-\e-V(x)+K_\al\ast\p\ast\mu_n(x)\le0,$$ 
i.e. $\tilde{u}$ is a subsolution to the equation
\begin{equation}\label{subsol bar la n}
-\Delta u+\frac{1}{\g}|\nabla u|^\g+\tilde{\la}-\e= V(x)-K_\al\ast\p\ast\mu_n(x).
\end{equation}
Since by definition (see \cite[Theorem 2.7 (i)]{CC})
$$\bar{\la}_n:=\sup\bigg\{\la\in\R\,\bigg|\,-\Delta u+\frac{1}{\g}|\nabla u|^\g+\la= V(x)-K_\al\ast\mu_n\ast\p(x)\quad\text{has a subsolution in } C^2(\R^N)\bigg\}$$
it must be $\bar{\la}_n\ge\tilde{\la}-\e$, which yields a contradiction.
Therefore $\tilde{\la}=\la^{(1)}$. By Proposition \ref{prop fix} $ii)$ $\tilde{u}$ is unique up to addition of constants, namely there exists $c\in\R$ such that $\tilde{u}=u^{(1)}+c$, it follows that $\nabla\tilde{u}=\nabla u^{(1)}$. Once we have the sequence of function $u_n$, we construct the sequence $(m_n)_n\subset \mathcal{F}(A_{\xi,M,C})$ such that for every $n\in\mathbb{N}$ fixed, it holds
$$\begin{cases}
-\Delta m_n-\mathrm{div}(m_n\,\nabla u_n|\nabla u_n|^{\g-2})=0\\
\int_{\R^N}m_n=M,\quad m_n\ge0
\end{cases}.$$
\noindent From Lemma \ref{pre comp}, up to extracting a subsequence
$$m_n\to m^{(1)}\,\,\,\text{in}\,\,\,L^{\bar{p}}(\R^N)$$
where $m^{(1)}\in A_{\xi,M,C}$. Since $\nabla u_n|\nabla u_n|^{\g-2}\to \nabla \tilde  u|\nabla \tilde u|^{\g-2}$ locally  uniformly in $\R^N$, we get that  $m^{(1)}$ is a weak solution to 
$$-\Delta m-\mathrm{div}(m\nabla\tilde{u}|\nabla\tilde{u}|^{\g-2})=0$$
that has $\tilde{m}=\mathcal{F}(\tilde{\mu})$ as unique solution. This proves that $m_n\to\tilde{m}$ in $L^{\bar{p}}(\R^N)$.
\endproof

We are ready to prove the following result on existence of solutions to the regularised MFG system \eqref{1 phi}.

\begin{theorem}\label{esist MFGreg}
We get the following results:
\begin{itemize}
    \item[i.] if $N-\g'<\al<N$ then, for every $M>0$ the MFG system \eqref{1 phi} admits a classical solution;
    \item[ii.]  if $N-2\g'\le\al\leq N-\g'$ then, there exists a positive real value $M_0=M_0(N,\al,\g,C_V,b)$ such that if $M\in(0, M_0)$ the MFG system \eqref{1 phi} admits a classical solution.
\end{itemize}
\end{theorem}

\proof
From Lemma \ref{A chiuso}, Lemma \ref{K in se}, Lemma \ref{pre comp} and Lemma \ref{F continua} assumptions of Theorem \ref{cor11.2GT} are verified, hence the map $\mathcal{F}$ has a fixed point $m_\p$. The fixed point $m_\p$ together with the couple $(u_\p,\bar{\la}_\p)=F(m_\p)$ obtained solving the Hamilton-Jacobi-Bellman equation with Riesz potential term equal to $K_\al\ast m_\p\ast\p$, will be a solution to the MFG system \eqref{1 phi}.
\endproof


\subsection{Limiting procedure}

Let $(\p_k)_k$ be a sequence of standard symmetric mollifiers approximating the unit as $k\to+\infty$. For every $k\in\mathbb{N}$ (under the additional assumption that the constraint mass $M$ is sufficiently small in the case when $N-2\g'<\al\le N-\g'$) from Theorem \ref{esist MFGreg} we can construct a classical solution $(u_k,m_k,\bar{\la}_k)$ to the corresponding regularised MFG system \eqref{1reg}. Our aim now is passing to the limit as $k\to+\infty$ and prove that $(u_k,m_k,\bar{\la}_k)$ converges to a solution of the MFG system \eqref{1}.

We need some preliminary apriori estimates. 

\begin{lemma}\label{lemma stime unif in k}
Let $\al\in(N-2\g',N)$ and $(u_k,m_k,\bar{\la}_k)$ be a solution to the regularized MFG system \eqref{1reg} as constructed in  Theorem \ref{esist MFGreg}. Then, there exist $C_1,C_2, C_3$ positive constants independent of $k$ such that  
\[\|m_k\|_{L^\infty(\R^N)}\le C_1, \quad\forall k\in \mathbb{N}\]
\[|\bar{\la}_k|\le C_2 \] 
and 
\begin{equation}\label{stime unif limiting}
|\nabla u_k|\le C_3(1+|x|^\frac{b}{\g})\qquad \qquad|\Delta u_k|\le C_3(1+|x|^b).
\end{equation}
\end{lemma}

\proof 
Note that  if $m\in L^1(\R^N)\cap L^\infty(\R^N)$ by Theorem \ref{holderRiesz}, we have that 
\begin{align*}
\|K_\al\ast\p_k\ast m\|_{L^\infty(\R^N)}&\le C_{N,\al,r,s}(\|\p_k\ast m\|_{L^r(\R^N)} +\|\p_k\ast m\|_{L^s(\R^N)})\\
&\le C_1\|m\|_{L^r(\R^N)} +C_2\|m\|_{L^s(\R^N)}
\end{align*}
for every $r\in\left(\frac{N}{\al},+\infty\right]$ and $s\in\left[1,\frac{N}{\al}\right)$. So, since $m_k\in L^1(\R^N)\cap L^\infty(\R^N)$ with $\|m_k\|_{\frac{2N}{N+\alpha}}\leq \xi$,   and $u_k$ are bounded from below, we may apply Theorem \ref{teo4.1new} with $W_k\equiv V$, $G_{k,\al}[m]=K_\al\ast\p_k\ast m$ and $q=\frac{2N}{N+\alpha}>\frac{N}{\alpha+\gamma'}$ and conclude the uniform $L^\infty$ bounds on $m_k$.\\
Now, by  Proposition \ref{prop fix}, we get that $\bar{\la}_k$ are equibounded in $k$ and that 
$$|\nabla u_k(x)|\le C(1+|x|^\frac{b}{\g})\qquad |\Delta u_k|\le C(1+|x|^b)$$ where $C$ is independent of $k$.  
\endproof

\proof[Proof of Theorem \ref{teo esist fix}]
Since for any $k\in\mathbb{N}$, $u_k$ is a classical solution to the HJB equation
$$-\Delta u+\frac{1}{\g}|\nabla u|^\g+\bar{\la}_k=V(x)-K_\al\ast\p_k\ast m_k$$ by Lemma \ref{lemma stime unif in k}  and elliptic estimates (refer to \cite[Theorem 8.32]{GT}) applied to $v_k(x):=u_k(x)-u_k(0)$, we obtain that for every $K\subset \subset \R^N$ and $\theta\in(0,1]$
$$\|v_k\|_{C^{1,\theta}_{loc}(K)}\le C\qquad\text{uniformly with respect to }k$$
hence up to extracting a subsequence 
$$v_k\to \bar{u} \qquad\text{locally uniformly in $C^1$ on compact sets}.$$
Similarly, since $m_k$ weak solution to $-\Delta m-\mathrm{div}(m|\nabla u_k|^{\g-2}\nabla u_k)=0$, for every $\phi\in C^\infty_0(K)$ it holds
$$\bigg|\int_K m_k\Delta\phi\,dx\bigg|\le\|\nabla\phi\|_{L^1(K)}\|m_k|\nabla u_k|^{\g-1}\|_{L^\infty(K)}.$$
Using  the uniform $L^\infty$ estimates on $m_k$ and  the estimates \eqref{stime unif limiting},  by Proposition \ref{prop2.4CC} and Sobolev embedding, we get that for every $\beta\in(0,1)$
$$\|m_k\|_{C^{0,\beta}(K)}\le C \qquad\text{uniformly with respect to }k$$
so up to extracting a subsequence 
$$m_k\to \bar{m} \qquad\text{locally uniformly}.$$
Since the values of $\bar{\la}_k$ are equibounded with respect to $k$, we have that $\bar{\la}_k\to \bar{\la}$ up to a subsequence.
Again recalling that $\int V(x)m_k\leq C$ uniformly in $k$, we conclude by  Vitali Convergence Theorem that $m_k\to\bar{m}$ in $L^1(\R^N)$ and hence $\int_{\R^N}\bar{m}=M$. From the strong convergence in $L^1(\R^N)$ and the uniform $L^\infty$ estimates, we obtain also that 
$$m_k\to\bar{m} \quad\text{in }L^p(\R^N)$$
for every $p\in[1,+\infty)$. We finally have that 
$$K_\al\ast\p_k\ast m_k\to K_\al\ast m \quad \text{locally uniformly}.$$
We can pass to the limit and obtain that $(\bar{u},\bar{m},\bar{\la})$ is a solution to the MFG system \eqref{1}.
\endproof


\end{document}